\documentclass[11pt]{article}
\usepackage{url,epsf,amsmath,amssymb,amsfonts,graphicx,tikz}
\usepackage{color}
\usepackage{hyperref}
\usepackage[title]{appendix}

\newtheorem{theorem}{\bf Theorem}[section]
\newtheorem{corollary}[theorem]{\bf Corollary}

\newtheorem{lemma}[theorem]{\bf Lemma}

\newtheorem{remark}[theorem]{\bf Remark}
\newtheorem{problem}[theorem]{\bf Problem}

\newcommand{\proof}{\noindent{\bf Proof.\ }}
\newcommand{\qed}{\hfill $\Box$ \bigskip}
\newcommand{\cp}{\,\square\,}

\newcommand{\gp}{{\rm gp}}

\newcommand{\n}{\#{\rm gp}}
\newcommand{\pr}{p_{_{P_r}}}
\newcommand{\ps}{p_{_{P_s}}}

\begin{document}

\title{On general position sets in Cartesian products}

\author{
Sandi Klav\v zar $^{a,b,c}$
\and
Bal\'azs Patk\'os $^d$
\and
Gregor Rus $^{c,e}$
\and
Ismael G. Yero $^f$
}

\date{}

\maketitle

\begin{center}
$^a$ Faculty of Mathematics and Physics, University of Ljubljana, Slovenia \\
\medskip

$^b$ Faculty of Natural Sciences and Mathematics, University of Maribor, Slovenia \\
\medskip

$^{c}$ Institute of Mathematics, Physics and Mechanics, Ljubljana, Slovenia \\
\medskip

$^d$ Alfr\'ed R\'enyi Institute of Mathematics, Hungarian Academy of Sciences, Budapest

\medskip

$^e$ Faculty of Organizational Sciences, University of Maribor, Slovenia\\
\medskip

$^{f}$ Departamento de Matem\'aticas, Universidad de C\'adiz, Algeciras, Spain
\end{center}

\begin{abstract}
The general position number ${\rm gp}(G)$ of a connected graph $G$ is the cardinality of a largest set $S$ of vertices such that no three distinct vertices from $S$ lie on a common geodesic; such sets are refereed to as gp-sets of $G$.  The general position number of cylinders $P_r\cp C_s$ is deduced. It is proved that ${\rm gp}(C_r\cp C_s)\in \{6,7\}$ whenever $r\ge s \ge 3$, $s\ne 4$, and $r\ge 6$.  A probabilistic lower bound on the general position number of Cartesian graph powers is achieved. Along the way a formula for the number of gp-sets in $P_r\cp P_s$, where $r,s\ge 2$, is also determined.
\end{abstract}

\noindent
{\bf E-mails}: sandi.klavzar@fmf.uni-lj.si, patkos.balazs@renyi.hu, gregor.rus4@um.si, \\ ismael.gonzalez@uca.es


\medskip\noindent
{\bf AMS Subj. Class.}: 05C12, 05C76, 05C30, 05D40

\section{Introduction}
\label{sec:intro}

Points in general position in the plane are basic objects in classical geometry as well as in computational geometry. The concept naturally extends to arbitrary metric spaces and asking for the maximum number of points in general position appears as a natural problem. The problem to find the maximum number of points that can be placed in the $n \times n$ grid so that no three points lie on a line is known as the {\em no-three-in-line problem} and was posed back in 1917 by  Dudeney~\cite{dudeney-1917}. After a century the problem remains open, cf.~\cite{ku-2018, misiak-2016, PoWo07}. Closely related problems are investigated in discrete geometry, see~\cite{cardinal-2017, PaWo13}.

Recently, the general position problem has been studied for metric spaces generated by graphs. The problem was independently introduced in~\cite{manuel-2018a, ullas-2016}, the present terminology and formalism are from~\cite{manuel-2018a}. If $G = (V(G), E(G))$ is a graph, then $S\subseteq V(G)$ is a {\em general position set} if $d_G(u,v) \ne d_G(u,w) + d_G(w,v)$ holds for every $\{u,v,w\}\in \binom{S}{3}$, where $d_G(x,y)$ denotes the shortest-path distance in $G$, and $\binom{S}{3}$ the set of all $3$-subsets of $S$. Equivalently, no three vertices lie on a common geodesic. We also say that the vertices from $S$ {\em lie in general position}. The {\em general position problem} is to find a largest general position set of $G$, the order of such a set is the {\em general position number} $\gp(G)$ of $G$. A general position set of $G$ of order $\gp(G)$ will be called a {\em gp-set}.

Following the seminal papers, the general position problem has been investigated in a sequence of papers~\cite{anand-2019+, ghorbani-2019, klavzar-2019+, manuel-2018b, neethu-2021, patkos-2020, thomas-2020, yao-2020+}. As it happens, in the special case of hypercubes, the general position problem was studied back in 1995 by K\"orner~\cite{korner-1995} related to some coding theory problems.  In this paper, asymptotic lower and upper bounds were proved on the gp-number of hypercubes, and several closely related problems (cf.\ Section~\ref{sec:3D-and-hypercubes}) were considered. The lower bound from~\cite{korner-1995} was improved in~\cite{randriambololona-2013}.

The results from~\cite{manuel-2018b} on the general position problem in interconnection networks with the emphasis on grid graphs were a starting motivation for the present study. One of the main results of~\cite{manuel-2018b} asserts that if $P_\infty$ denotes the two-way infinite path, then $\gp(P_\infty\cp P_\infty) = 4$, and consequently $\gp(P_r\cp P_s) = 4$ for $r,s\ge 3$. The non-trivial part of this result (that $\gp(P_r\cp P_s) \le 4$ holds) was proved using the so-called Monotone Geodesic Lemma which was in turn derived from the celebrated Erd\"os-Szekeres theorem, cf.~\cite[Theorem 1.1]{bukh-2014}. The result $\gp(P_\infty\cp P_\infty) = 4$ was recently generalized in two directions. In~\cite{klavzar-2021}, it is proved that $\gp(P_\infty^{\cp, k}) = 2^{2^{k-1}}$, while in~\cite{tian-2020+}, it is demonstrated that if $T_1$ and $T_2$ are arbitrary trees, then $\gp(T_1\cp T_2) = \gp(T_1) + \gp(T_2)$. The general position number of several other Cartesian product graphs is studied in~\cite{tian-2021}.

In the rest of this section we prepare material needed later on. In the next section we have a closer look at the structure of gp-sets in grid graphs. We need to do it for the proof of the subsequent theorem on cylinder graphs (Theorem~\ref{thm:cylinders}). As a side result we determine the number of gp-sets in $P_r\cp P_s$ for every $r,s\ge 2$, a result that could be of independent interest. In Section~\ref{sec:cyclenders} we then determine $\gp(P_r\cp C_s)$ for every $r\ge 2$ and $s\ge 3$. In the subsequent section we prove that if $3\le s\le r$, $4 \ne s$, and $r\ge 6$, then $\gp(C_r\cp C_s) \in  \{6,7\}$. We also discuss exact values and in particular prove that $\gp(C_3\cp C_s) = 6$ holds for $s\ge 6$. Motivated by the results of~\cite{korner-1995}, we consider in Section~\ref{sec:3D-and-hypercubes} how to apply the probabilistic method to obtain asymptotic lower bounds on the gp-number of Cartesian powers of graphs.

\subsection{Preliminaries}

For a positive integer $k$ we will use the notation $[k] = \{1,\ldots, k\}$ and $[k]_0 = \{0,\ldots k-1\}$. The {\em Cartesian product} $G\cp H$ of graphs $G$ and $H$ has the vertex set $V(G\cp H) = V(G)\times V(H)$, vertices $(g,h)$ and $(g',h')$ are adjacent if either $gg'\in E(G)$ and $h=h'$, or $g=g'$ and $hh'\in E(H)$. Given a vertex $h\in V(H)$, the subgraph of $G\cp H$ induced by the set $\{(g,h):\ g\in V(G)\}$, is a {\em $G$-layer} and is denoted by $G^h$. $H$-layers $^gH$ are defined analogously. Each $G$-layer and $H$-layer is  isomorphic to $G$ and $H$, respectively. If $X\subseteq V(G\cp H)$, then the {\em projection} $p_G(X)$ of $X$ to $G$ is the set $\{g\in V(G):\ (g,h)\in X\ {\rm for\ some}\ h\in V(H)\}$. The projection $p_H(X)$ of $X$ to $H$ is defined analogously. The $k$-tuple Cartesian product of a graph $G$ by itself, alias Cartesian power of $G$, is denoted by $G^{\cp, k}$. This is well-defined since the Cartesian product operation is associative. For more on the Cartesian product see~\cite{imrich-2008}. As stated above, the following result was the primary motivation for the present paper.

\begin{theorem} {\rm \cite{manuel-2018b}}
	\label{thm:grids}
	If $r\ge 3$ and $s\ge 3$, then $\gp(P_r\cp P_s) = 4$.
\end{theorem}

A subgraph $H$ of a graph $G$ is {\em isometric} if $d_H(u,v) = d_G(u,v)$ holds for all $u,v\in V(H)$. A set of subgraphs $\{H_1,\ldots, H_k\}$ of a graph $G$ is an {\em isometric cover} of $G$ if each $H_i$, $i\in [k]$, is isometric in $G$ and $\bigcup_{i=1}^k V(H_i) = V(G)$.

\begin{theorem} {\rm \cite[Theorem 3.1]{manuel-2018a}}
	\label{thm:isometric-cover}
	If $\{H_1,\ldots, H_k\}$ is an isometric cover of $G$, then
	$$\gp(G) \le \sum_{i=1}^k \gp(H_i)\,.$$
\end{theorem}

If $G$ is a connected graph, $S\subseteq V(G)$, and ${\cal P} = \{S_1, \ldots, S_p\}$ is a partition of $S$, then ${\cal P}$ is \emph{distance-constant} (alias ``distance-regular''~ \cite[p.~331]{kante-2017}) if for any $i,j\in [p]$, $i\ne j$, the distance $d_G(u,v)$, where $u\in S_i$ and $v\in S_j$, is independent of the selection of $u$ and $v$. This distance is then the distance $d_G(S_i,S_j)$ between the parts $S_i$ and $S_j$. A distance-constant partition ${\cal P}$ is {\em intransitive} if $d_G(S_i, S_k) \ne d_G(S_i, S_j) + d_G(S_j,S_k)$ holds for distinct indices $i,j,k\in [p]$.

\begin{theorem} {\rm \cite[Theorem 3.1]{anand-2019+}}
	\label{thm:gpsets}
	Let $G$ be a connected graph. Then $S\subseteq V(G)$ is a general position set if and only if the components of the subgraph induced by $S$ are complete subgraphs, the vertices of which form an intransitive, distance-constant partition of $S$.
\end{theorem}

Suppose that $G$ is a connected bipartite graph and a general position set $S$ contains two adjacent vertices $x$ and $y$. Then Theorem~\ref{thm:gpsets} implies that $|S| = 2$, because no other vertex of $G$ can be at the same distance to $x$ and $y$. We state this observation for later use.

\begin{corollary}
	\label{cor:bipartite}
	If $G$ is a bipartite graph with $\gp(G) \ge 3$, then every gp-set of $G$ is an independent set.
\end{corollary}

\section{Enumeration of gp-sets in grids}
\label{sec:grids}

In this section we take a closer look at the structure of gp-sets in grids. We have two reasons to do it, the first is that this insight will be used in the proof of Theorem~\ref{thm:cylinders}, the second is that we are able to enumerate the gp-sets in grids as follows.

\begin{theorem}
	\label{thm:grids-enumerate}
	If $2\le r\le s$, then
	$$
	\n(P_r\cp P_s) = \left\{
	\begin{array}{ll}
	\vspace*{2mm}
	6; & r = s = 2\,, \\
	\vspace*{2mm}
	\displaystyle{\frac{s(s-1)(s-2)}{3}}; & r = 2, s\ge 3\,, \\
	\displaystyle{\frac{rs(r - 1)(r - 2)(s - 1)(s - 2)(r(s - 3) - s + 7)}{144}}; & r, s\ge 3\,.
	\end{array}
	\right.
	$$
\end{theorem}

\proof
Set $V(P_n) = [n]$. If $r = s = 2$, then the assertion is clear since $P_2\cp P_2 = C_4$. Let next $r=2$ and $s\ge 3$. It is straightforward to see that $\gp(P_2\cp P_s) = 3$. Moreover, if $X$ is a gp-set of $P_2\cp P_s$, then $X$ has one vertex in one of the two $P_s$-layers and two vertices in the other $P_s$-layer, say $X = \{(1,i), (2,j), (2,k)\}$, where $j<k$. Since $X$ is a gp-set we infer that $j < i < k$. Hence, for a given vertex $(1,i)$, there are $i-1$ possibilities to select the vertex $(2,j)$, and $s-i$ possibilities for the vertex $(2,k)$. The same holds if $X$ has two vertices in $^1P_s$ and one vertex in $^2P_s$. From this it follows that
$$\n(P_2\cp P_s) = 2\cdot \sum_{i=1}^s (i-1)(s-i) = \frac{s(s-1)(s-2)}{3}\,.$$
Suppose in the rest that $r,s\ge 3$, so that $\gp(P_r\cp P_s) = 4$ by Theorem~\ref{thm:grids}. Hence by Corollary~\ref{cor:bipartite}, every gp-set is an independent set (of cardinality $4$). Let $X$ be an arbitrary such set and assume first that $|\ps(X)| = 2$. Then, clearly, $X$ has two vertices in one $P_r$-layer and two vertices in another $P_r$-layer. Let $(i,j)\in X$ be a vertex that has the smallest first coordinate among the vertices of $X$. Then $(i,j)$ and the two vertices of $X$ from the $P_r$-layer not containing $(i,j)$ lie on a common geodesic. Analogously, $X$ cannot be a general position set if $|\pr(X)| = 2$. Since also $|\pr(X)| = 1$ or $|\ps(X)| = 1$ are not possible, we only need to distinguish the following two cases.

\medskip\noindent
{\bf Case 1}: $|\pr(X)| = 4$ and $|\ps(X)| = 4$. \\
Let $\pr(X) = \{a, b, c, d\}$, where $a < b < c < d$, and let $\ps(X) = \{a', b', c',  d'\}$,  where $a' < b' < c' < d'$. Then in the set $\pr(X) \times \ps(X)$ there are $4!$ different $4$-sets of vertices that project onto both $\pr(X)$ and $\ps(X)$. They can be described with permutations $\pi$ of $\ps(X)$. That is, if $\pi: \ps(X) \rightarrow \ps(X)$ is a bijection, then the corresponding gp-set of vertices of $P_r\cp P_s$ is $S_\pi = \{(a, \pi(a')), (b, \pi(b')), (c, \pi(c')), (d, \pi(d'))\}$. Now, by the metric structure of $P_r\cp P_s$ (cf.~\cite{manuel-2018b}), $S_\pi$ is a general position set if and only if the sequence $(\pi(a'), \pi(b'), \pi(c'), \pi(d'))$ contains no monotone subsequence of length $3$. By a direct inspection we find that if $\pi(a') = a'$ or if $\pi(a') = d'$, then we get no general position sets. If $\pi(a') = b'$, then exactly the sequences $(b', a', d', c')$ and $(b', d', a', c')$ yield general position sets. Symmetrically, if $\pi(a') = c'$, then exactly the sequences $(c', a', d', b')$ and $(c', d', a', b')$ yield general position sets. Hence, if $|\pr(X)| = 4$ and $|\ps(X)| = 4$, then there are exactly $4\binom{r}{4}\binom{s}{4}$ gp-sets.

\medskip\noindent
{\bf Case 2}: $|\pr(X)| = 3$ (and $|\ps(X)| = 3$ or $|\ps(X)| = 4$). \\
Let $\pr(X) = \{a, b, c\}$, where two vertices from $X$ project to $a$, say $(a,a'), (a,b')\in X$, where $a' < b'$. Let $(x,x')$ be a vertex of $P_r\cp P_s$, where $x'\le a'$ and $(x,x')\ne (a,a')$. Then $d((x,x'), (a,b')) = d((x,x'), (a,a')) + d((a,a'), (a,b'))$, which means that $(x,x')\notin X$. Similarly, if $(x,x')$ is a vertex of $P_r\cp P_s$ with $x'\ge b'$ and $(x,x')\ne (a,b')$, then also  $(x,x')\notin X$. We have thus shown that
$$X\cap ([r] \times \{1,\ldots, a'\}) = \{(a,a')\}\ {\rm and}\ X\cap ([r] \times \{b',\ldots, s\}) = \{(a,b')\}\,.$$
Let $X = \{(a,a'), (a,b'), (x,x'), (y,y')\}$, where $x\ne y$. By a similar argument as above we see that, without loss of generality,
\begin{align*}
(x,x') & \in \{1,\ldots, a-1\}\times \{a'+1,\ldots, b'-1\}\ {\rm and}\\
(y,y') & \in \{a+1,\ldots, r\}\times \{a'+1,\ldots, b'-1\}\,.
\end{align*}
Since the vertices $(x,x')$ and $(y,y')$ are arbitrary vertices from $\{1,\ldots, a-1\}\times \{a'+1,\ldots, b'-1\}$ and $\{a+1,\ldots, r\}\times \{a'+1,\ldots, b'-1\}$, respectively, for fixed $a, a', b'$ we obtain precisely
$$[(b'-a'-1)(a-1)]\cdot [(b'-a'-1)(r-a)] = (b'-a'-1)^2(a-1)(r-a)$$
gp-sets. To obtain the number of all gp-sets in this case, we need to sum up over all coordinates $a$ of the factor $P_r$, onto which two vertices from $X$ project. Consequently, the number of gp-sets in Case 2 is
$$\sum_{a=1}^r \sum_{a'=1}^s \sum_{b'=a'+1}^s [(b'-a'-1)^2(a-1)(r-a)] = \frac{rs(r^2 - 3r + 2)(s^3-4s^2 + 5s - 2)}{72}\,.$$

\noindent
By the above two cases, if $r, s\ge 3$, then
\begin{align*}
\n(P_r\cp P_s) & = 4\binom{r}{4}\binom{s}{4} + \frac{rs(r^2 - 3r + 2)(s^3-4s^2 + 5s - 2)}{72} \\
& = \frac{rs(r - 1)(r - 2)(s - 1)(s - 2)(r(s - 3) - s + 7)}{144}\,,
\end{align*}
which is the claimed expression.
\qed

If $r=3$ and $s\ge 3$, then Theorem~\ref{thm:grids-enumerate} yields
$$\n(P_3\cp P_s) = \frac{s(s - 2)(s - 1)^2}{12}\,,$$
which, after substituting $s$ with $s+1$ gives the sequence A002415 from OEIS~\cite{sloane}. In addition, the case $r=2$ and $s\ge 3$ yields the sequence A007290.

\section{Cylinders}
\label{sec:cyclenders}

In this section we determine the general position number of cylinders. For this task, the following function will be useful. If $G$ is a connected graph and $X\subseteq V(G)$ is a general position set, then
$$F(X) = \{u\in V(G)-X:\ X\cup \{u\}\ {\rm is\ not\ a\ general\ position\ set}\}\,.$$
If $X=\{x,y\}$, we will simplify the notation $F(\{x,y\})$ to $F(x,y)$.

Set $V(P_r)=[r]_0$ and $V(C_s)=[s]_0$. From now on, operations with the integers in $V(C_s)$ are done modulo $s$.

\begin{lemma}
	\label{lem:general}
	Let $r\ge 2$, $s\ge 3$, and let $S$ be a general position set of the cylinder graph $P_r\cp C_s$. Then the following assertions hold.
	\begin{itemize}
		\item[(i)] If $|S|\ge 5$, then $S$ is an independent set.
		\item[(ii)] If $|S| \ge 4$, then $|S\cap V(^iC_s)|\le 2$ for every $i\in [r]_0$.
		\item[(iii)] If $|S\cap V(^iC_s)| = 2$ for some $i\in [r]_0$, then $|S|\le 4$.
		\item[(iv)] If $r\ge 6$, $|S| = 5$, and $|S\cap V(^iC_s)| \le 1$ for every  $i\in [r]_0$, then $\gp(P_5\cp C_s) \ge 5$.
	\end{itemize}
\end{lemma}

\proof
(i) Suppose $|S|\ge 5$ and $S$ is not independent. If $(i,k), (i+1,k)\in S$, then we observe that $F((i,k), (i+1,k))=V(P_r\cp C_s)-\{(i,k), (i+1,k)\}$, which means that $S = \{(i,k), (i+1,k)\}$, a contradiction. On the other hand, if $(i,k), (i,k+1)\in S$, then either $F((i,k), (i,k+1))=V(P_r\cp C_s)-\{(i,k), (i,k+1)\}$ (when $s$ is even), or $F((i,k), (i,k+1))=V(P_r\cp C_s)-(\{(i,k), (i,k+1)\}\cup (V(P_r)\times \{j\}))$ (when $s$ is odd), where $j$ is the vertex which is on $C_s$ opposite to $k$ and $k+1$. The first possibility directly leads to a contradiction. For the second one, since every $P_r$-layer, being an isometric subgraph, contributes at most two vertices to a general position set of $P_r\cp C_s$, it follows that $|S|\le 4$, which is again a contradiction. Consequently $S$ must be an independent set.

(ii) The result follows directly from the following fact. If $\{(i,k_1),(i,k_2),(i,k_3)\} \subseteq S$, then (since $k_1$, $k_2$, and $k_3$ are distinct) $F(\{(i,k_1),(i,k_2),(i,k_3)\})=V(P_r\cp C_s)-\{(i,k_1),(i,k_2),(i,k_3)\}$, which means that $|S| = 3$.

(iii) Let $i\in [r]_0$ be such that $|S\cap V(^iC_s)| = 2$. We may assume without loss of generality that $S\cap V(^iC_s) = \{(i,0), (i,j)\}$, where $j\le \lfloor s/2 \rfloor$. Then
$$F((i,0), (i,j)) = [r]_0 \times (\{j-\lfloor s/2\rfloor,\dots, 0\}\cup \{j,\dots, \lfloor s/2\rfloor\})\cup (\{i\}\times I)-\{(i, 0), (i, j)\},$$
where $I=[s]_0$ if $s$ is even and $j=s/2$; or $I=\{0,\dots, j\}$ otherwise.

If $j=1$, then $F((i,0), (i,j))$ equals $V(P_r\cp C_s) - \{(i, 0), (i, j)\}$ when $s$ is even, or $F((i,0), (i,j))$ equals $V(P_r\cp C_s)-(([r]_0\times \{\lceil s/2\rceil\}) \cup \{(i, 0), (i, j)\})$ when $s$ is odd. In the first situation we clearly have $|S|\le 2$. In the latter one, the set of vertices $[r]_0\times \{\lceil s/2\rceil\}$ could contain at most two vertices of $S$, since it induces an isometric path in $P_r\cp C_s$, and so $|S|\le 4$. We may assume in the rest that $j\ge 2$. Note that $V(P_r\cp C_s)-(F((i,0), (i,j))\cup \{(i,0), (i,j)\})=Y_1\cup Y_2\cup Y_3\cup Y_4$, where
\begin{align*}
  Y_1= & \,[i]_0 \times \{1,\dots, j-1\} \\
  Y_2= & \,[i]_0 \times \{\lfloor s/2\rfloor+1,\dots, j-\lfloor s/2\rfloor-1\} \\
  Y_3= & \,([r]_0-[i+1]_0) \times \{1,\dots, j-1\} \\
  Y_4= & \,([r]_0-[i+1]_0) \times \{\lfloor s/2\rfloor+1,\dots, j-\lfloor s/2\rfloor-1\}.
\end{align*}

Consider two vertices $(i',j'),(i'', j'')\in S-\{(i,0), (i,j)\}$. (If there are no two such vertices, then clearly $|S|\le 3$.) If either $(i',j'),(i'', j'')\in Y_1$; or $(i',j'),(i'', j'')\in Y_2$; or $(i',j'),(i'', j'')\in Y_3$; or $(i',j'),(i'', j'')\in Y_4$, then it happens that $(i,0)\in F((i',j'), (i'',j''))$ or $(i,j)\in F((i',j'), (i'',j''))$, which is not possible. As a consequence, each of the sets $Y_1$, $Y_2$, $Y_3$, and $Y_4$ must contain at most one vertex of $S$.

Assume now there is a vertex $(i',j')\in Y_1\cap S$. Then the set $F((i',j'), (i,0))\cup F((i',j'), (i,j))$ contains the whole set $Y_4$. By symmetry, if there is a vertex $(i',j')\in Y_4\cap S$, then $Y_1\cap S = \emptyset$. Hence, $|(Y_1\cup Y_4) \cap S| \le 1$. Similarly we obtain that $|(Y_2\cup Y_3) \cap S| \le 1$. Therefore, the sets $Y_1$, $Y_2$, $Y_3$, and $Y_4$ can contain at most two vertices of $S$, which gives $|S|\le 4$.

(iv) Let $S = \{(i_k, j_k):\ k\in [5]_0\}$. Since  $|S\cap V(^iC_s)| \le 1$, the coordinates $i_k$ are pairwise different, hence we may assume without loss of generality that $i_0 < i_1 < i_2 < i_3 < i_4$. Set $S' = \{(k,j_k):\ k\in [5]_0\}$. We claim that $S'$ is a general position set of $G_5 = P_5\cp C_s$. Assume on the contrary that
$$d_{G_5}((p,j_p), (r,j_{r})) = d_{G_5}((p,j_p), (q,j_{q})) + d_{G_5}((q,j_q), (r,j_{r}))\,$$
for some $p,q,r\in [5]_0$, $p < q < r$. Since the distance function in Cartesian products is additive, we get that
$$d_{P_5}(p,r) + d_{C_s}(j_p,j_{r}) = d_{P_5}(p,q) + d_{C_s}(j_p,j_{q}) + d_{P_5}(q,r) + d_{C_s}(j_q,j_{r})\,.$$
Since $d_{P_5}(p,r) = d_{P_5}(p,q) + d_{P_5}(q,r) $, we thus have
$$d_{C_s}(j_p,j_{r}) = d_{C_s}(j_p,j_{q}) + d_{C_s}(j_q,j_{r})\,.$$
From this we get that in $G_r = P_r\cp C_s$,
\begin{align*}
d_{G_r}((i_p,j_p), (i_r,j_{r})) & = d_{P_r}(i_p,i_r) + d_{C_s}(j_p, j_{r}) \\
& = [d_{P_r}(i_p,i_q) + d_{P_r}(i_q,i_r)] + [d_{C_s}(j_p,j_{q}) + d_{C_s}(j_q,j_{r})] \\
& = [d_{P_r}(i_p,i_q) + d_{C_s}(j_p,j_{q})] + [d_{P_r}(i_q,i_r) + d_{C_s}(j_q,j_{r})] \\
& = d_{G_r}((i_p,j_p), (i_q,j_{q})) + d_{G_r}((i_q,j_q), (i_r,j_{r}))\,.
\end{align*}
This contradiction proves that $S'$ is a general position set of $P_5\cp C_s$. We conclude that $\gp(P_5\cp C_s) \ge 5$.
\qed

Note that Lemma~\ref{lem:general}(iv) allows us to map a general position set of cardinality $5$ in long cylinders to a general position set of the same cardinality in cylinders over $P_5$.

\begin{theorem}
	\label{thm:cylinders}
	If $r\ge 2$ and $s\ge 3$, then
	$$
	\gp(P_r\cp C_s) = \left\{
	\begin{array}{ll}
	3; & r = 2, s = 3\,, \\
	5; & r \ge 5,\ {\rm and}\ s=7\ {\rm or}\ s\ge 9\,, \\
	4; & {\rm otherwise}\,.
	\end{array}
	\right.
	$$
\end{theorem}

\proof
First, it is easy to verify that $\gp(P_2\cp C_3) = 3$.

Assume next that $r\le 4$ and suppose that there exists a general position set $S$ with $|S|\ge 5$. If $r=2$, this is not possible by Lemma~\ref{lem:general}(ii). Let next $r\in \{3,4\}$. Then there exists a $C_s$-layer $^iC_s$ with $|V(^iC_s)\cap S| \ge 2$. The case $|V(^iC_s)\cap S| > 2$ is not possible by Lemma~\ref{lem:general}(ii), while the case $|V(^iC_s)\cap S| = 2$ is excluded by Lemma~\ref{lem:general}(iii). Hence $\gp(P_r\cp C_s) \le 4$ for $r\in \{2,3,4\}$. It is straightforward to see that the set $\{(0,0), (1,1), (0, \lfloor s/2\rfloor), (1, \lfloor s/2\rfloor +1)\}$ is a general position set of
$P_r\cp C_s$ for $r\ge 2$ and $s\ge 4$. Moreover, if $s=3$, then the set $\{(0,1), (1,0), (1,2), (2,1)\}$ is a general position set of $P_r\cp C_3$ for $r\ge 3$. Hence $\gp(P_r\cp C_s) \ge 4$ for $r\in \{2,3,4\}$ and so $\gp(P_r\cp C_s) = 4$ for $r\in \{2,3,4\}$.

The general position set $\{(0,0), (1,2), (2,4), (3,6), (4,1)\}$ of $P_5\cp C_7$ demonstrates that $\gp(P_5\cp C_7) \ge 5$.

Suppose next that for some $r\ge 6$ the cylinder $P_r\cp C_8$ contains a general position set $S$ with $|S| = 5$. From Lemma~\ref{lem:general}(iii) it follows that $|S\cap V(^iC_8)| \le 1$ for every  $i\in [r]_0$. Hence the assumptions of Lemma~\ref{lem:general}(iv) are fulfilled, and we can map the general position set $S$ of cardinality $5$ in the cylinder $P_r\cp C_8$, $r\ge 6$, to a general position set of the same cardinality ($5$) in the cylinder over $P_5\cp C_8$, which indeed implies that $\gp(P_5\cp C_8) \ge 5$. However, we have checked by computer that $\gp(P_5\cp C_8) = 4$. Thus, we have obtained a contradiction. Therefore, $\gp(P_r\cp C_8) \le 4$ for $r\ge 5$ (the case $r=5$ was computationally made). Since clearly $\gp(P_r\cp C_8) \ge 4$, we conclude that $\gp(P_r\cp C_8) = 4$ for $r\ge 5$.

Suppose now that $r = 5$, $s\ge 9$, and consider the set
$$S = \{u_0=(0,1), u_1=(1,4), u_2=(2, \lfloor s/2\rfloor+2), u_3=(3,0), u_4=(4, 3)\}\,.$$
We claim that $S$ is a general position set. Note first that the vertices $u_0, u_1, u_3, u_4$ lie in a general position. Further,
\begin{align*}
d(u_2,u_0) & = (s - \lfloor s/2\rfloor - 1) + 2 = s - \lfloor s/2\rfloor + 1,\\
d(u_2,u_1) & = \lfloor s/2\rfloor - 1, \\
d(u_2,u_3) & = s - \lfloor s/2\rfloor - 1, \\
d(u_2,u_4) & = \lfloor s/2\rfloor + 1.
\end{align*}
Then $d(u_0,u_2) = s - \lfloor s/2\rfloor + 1 < \lfloor s/2\rfloor + 3 =
d(u_0,u_1) + d(u_1,u_2)$. Similarly we see that $u_2$ is not on a geodesic containing three vertices of $S$. Hence, $S$ is a general position set and thus $\gp(P_5\cp C_s) \ge 5$ for $s\ge 9$.

Note finally that the general position set for $P_5\cp C_7$ and the general position set for $P_5\cp C_s$, $s\ge 9$, are also general position sets for $P_r\cp C_7$, $r\ge 6$, and for $P_r\cp C_s$, $r\ge 6$, respectively. We conclude that $\gp(P_r\cp C_s) \ge 5$ for $r\ge 5$ and $s\ge 7$, $s\ne 8$.

It remains to prove that the above constructed general position sets of cardinality $5$ are gp-sets. Hence let $S$ be a gp-set of $P_r\cp C_s$, where $|S| \ge 5$. Note that $S$ is an independent set, by Lemma~\ref{lem:general}(i). Then make a partition of $V(P_r\cp C_s)$ into two sets $A_1$ and $A_2$ inducing two grids that are isometric subgraphs of $P_r\cp C_s$. Without loss of generality, we may assume $A_1=V(P_r)\times [\lceil s/2\rceil]_0$ and $A_2=V(P_r)\times ([s]_0-[\lceil s/2\rceil]_0)$. Now, let $S_1=S\cap A_1$ and $S_2=S\cap A_2$. Since $\gp(P_r\cp C_s)\ge 5$, it follows $|S_1|\ge 3$ or $|S_2|\ge 3$. Moreover, since $A_1$ and $A_2$ induce isometric grid graphs, Theorem~\ref{thm:grids} implies that $|S_1|\le 4$ and $|S_2|\le 4$. To simplify the notation, we write $s_d=\lfloor s/2\rfloor$ Note that here we do not need to change floor to ceil. We consider the following cases.

\begin{figure}[ht!]
	\centering
	\begin{tikzpicture}[scale=.75, transform shape]
	\node [draw, shape=circle,scale=0.7] (a1) at  (0,0) {};
	\node [draw, shape=circle,scale=0.7] (a2) at  (0,1) {};
	\node [draw, shape=circle,scale=0.7] (a3) at  (0,2) {};
	\node [draw, shape=circle,scale=0.7] (a4) at  (0,3) {};
	\node [draw, shape=circle,scale=0.7] (a5) at  (0,4) {};
	\node [draw, shape=circle,scale=0.7] (a6) at  (0,5) {};
	\node [draw, shape=circle,scale=0.7] (a7) at  (0,6) {};
	
	\node [draw, shape=circle,scale=0.7] (b1) at  (1,0) {};
	\node [draw, shape=circle,scale=0.7] (b2) at  (1,1) {};
	\node [draw, shape=circle,scale=0.7] (b3) at  (1,2) {};
	\node [draw, shape=circle,scale=0.7,fill=black] (b4) at  (1,3) {};
	\node [draw, shape=circle,scale=0.7] (b5) at  (1,4) {};
	\node [draw, shape=circle,scale=0.7] (b6) at  (1,5) {};
	\node [draw, shape=circle,scale=0.7] (b7) at  (1,6) {};
	
	\node [draw, shape=circle,scale=0.7] (c1) at  (2,0) {};
	\node [draw, shape=circle,scale=0.7,fill=black] (c2) at  (2,1) {};
	\node [draw, shape=circle,scale=0.7] (c3) at  (2,2) {};
	\node [draw, shape=circle,scale=0.7] (c4) at  (2,3) {};
	\node [draw, shape=circle,scale=0.7] (c5) at  (2,4) {};
	\node [draw, shape=circle,scale=0.7] (c6) at  (2,5) {};
	\node [draw, shape=circle,scale=0.7] (c7) at  (2,6) {};
	
	\node [draw, shape=circle,scale=0.7] (d1) at  (3,0) {};
	\node [draw, shape=circle,scale=0.7] (d2) at  (3,1) {};
	\node [draw, shape=circle,scale=0.7] (d3) at  (3,2) {};
	\node [draw, shape=circle,scale=0.7] (d4) at  (3,3) {};
	\node [draw, shape=circle,scale=0.7] (d5) at  (3,4) {};
	\node [draw, shape=circle,scale=0.7,fill=black] (d6) at  (3,5) {};
	\node [draw, shape=circle,scale=0.7] (d7) at  (3,6) {};
	
	\node [draw, shape=circle,scale=0.7] (e1) at  (4,0) {};
	\node [draw, shape=circle,scale=0.7] (e2) at  (4,1) {};
	\node [draw, shape=circle,scale=0.7] (e3) at  (4,2) {};
	\node [draw, shape=circle,scale=0.7] (e4) at  (4,3) {};
	\node [draw, shape=circle,scale=0.7,fill=black] (e5) at  (4,4) {};
	\node [draw, shape=circle,scale=0.7] (e6) at  (4,5) {};
	\node [draw, shape=circle,scale=0.7] (e7) at  (4,6) {};
	
	\node [draw, shape=circle,scale=0.7] (f1) at  (5,0) {};
	\node [draw, shape=circle,scale=0.7] (f2) at  (5,1) {};
	\node [draw, shape=circle,scale=0.7] (f3) at  (5,2) {};
	\node [draw, shape=circle,scale=0.7] (f4) at  (5,3) {};
	\node [draw, shape=circle,scale=0.7] (f5) at  (5,4) {};
	\node [draw, shape=circle,scale=0.7] (f6) at  (5,5) {};
	\node [draw, shape=circle,scale=0.7] (f7) at  (5,6) {};
	
	\node [scale=1] at (1,-0.5) {$a$};
	\node [scale=1] at (2,-0.5) {$b$};
	\node [scale=1] at (3,-0.5) {$c$};
	\node [scale=1] at (4,-0.5) {$d$};
	
	\node [scale=1] at (-1,1) {$b'$};
	\node [scale=1] at (-1,3) {$a'$};
	\node [scale=1] at (-1,4) {$d'$};
	\node [scale=1] at (-1,5) {$c'$};
	
	\end{tikzpicture}
	\hspace*{1cm}
	\begin{tikzpicture}[scale=.75, transform shape]
	\node [draw, shape=circle,scale=0.7] (a1) at  (0,0) {};
	\node [draw, shape=circle,scale=0.7] (a2) at  (0,1) {};
	\node [draw, shape=circle,scale=0.7] (a3) at  (0,2) {};
	\node [draw, shape=circle,scale=0.7] (a4) at  (0,3) {};
	\node [draw, shape=circle,scale=0.7] (a5) at  (0,4) {};
	\node [draw, shape=circle,scale=0.7] (a6) at  (0,5) {};
	\node [draw, shape=circle,scale=0.7] (a7) at  (0,6) {};
	
	\node [draw, shape=circle,scale=0.7] (b1) at  (1,0) {};
	\node [draw, shape=circle,scale=0.7] (b2) at  (1,1) {};
	\node [draw, shape=circle,scale=0.7] (b3) at  (1,2) {};
	\node [draw, shape=circle,scale=0.7] (b4) at  (1,3) {};
	\node [draw, shape=circle,scale=0.7,fill=black] (b5) at  (1,4) {};
	\node [draw, shape=circle,scale=0.7] (b6) at  (1,5) {};
	\node [draw, shape=circle,scale=0.7] (b7) at  (1,6) {};
	
	\node [draw, shape=circle,scale=0.7] (c1) at  (2,0) {};
	\node [draw, shape=circle,scale=0.7] (c2) at  (2,1) {};
	\node [draw, shape=circle,scale=0.7] (c3) at  (2,2) {};
	\node [draw, shape=circle,scale=0.7] (c4) at  (2,3) {};
	\node [draw, shape=circle,scale=0.7] (c5) at  (2,4) {};
	\node [draw, shape=circle,scale=0.7,fill=black] (c6) at  (2,5) {};
	\node [draw, shape=circle,scale=0.7] (c7) at  (2,6) {};
	
	\node [draw, shape=circle,scale=0.7] (d1) at  (3,0) {};
	\node [draw, shape=circle,scale=0.7,fill=black] (d2) at  (3,1) {};
	\node [draw, shape=circle,scale=0.7] (d3) at  (3,2) {};
	\node [draw, shape=circle,scale=0.7] (d4) at  (3,3) {};
	\node [draw, shape=circle,scale=0.7] (d5) at  (3,4) {};
	\node [draw, shape=circle,scale=0.7] (d6) at  (3,5) {};
	\node [draw, shape=circle,scale=0.7] (d7) at  (3,6) {};
	
	\node [draw, shape=circle,scale=0.7] (e1) at  (4,0) {};
	\node [draw, shape=circle,scale=0.7] (e2) at  (4,1) {};
	\node [draw, shape=circle,scale=0.7,fill=black] (e3) at  (4,2) {};
	\node [draw, shape=circle,scale=0.7] (e4) at  (4,3) {};
	\node [draw, shape=circle,scale=0.7] (e5) at  (4,4) {};
	\node [draw, shape=circle,scale=0.7] (e6) at  (4,5) {};
	\node [draw, shape=circle,scale=0.7] (e7) at  (4,6) {};
	
	\node [draw, shape=circle,scale=0.7] (f1) at  (5,0) {};
	\node [draw, shape=circle,scale=0.7] (f2) at  (5,1) {};
	\node [draw, shape=circle,scale=0.7] (f3) at  (5,2) {};
	\node [draw, shape=circle,scale=0.7] (f4) at  (5,3) {};
	\node [draw, shape=circle,scale=0.7] (f5) at  (5,4) {};
	\node [draw, shape=circle,scale=0.7] (f6) at  (5,5) {};
	\node [draw, shape=circle,scale=0.7] (f7) at  (5,6) {};
	
	\node [scale=1] at (1,-0.5) {$a$};
	\node [scale=1] at (2,-0.5) {$b$};
	\node [scale=1] at (3,-0.5) {$c$};
	\node [scale=1] at (4,-0.5) {$d$};
	
	\node [scale=1] at (-1,1) {$c'$};
	\node [scale=1] at (-1,2) {$d'$};
	\node [scale=1] at (-1,4) {$a'$};
	\node [scale=1] at (-1,5) {$b'$};
	
	\end{tikzpicture}
	\caption{Two possible configurations of the set $S_1$ (edges of the grid have not been drawn).}
	\label{fig-S_1}
\end{figure}
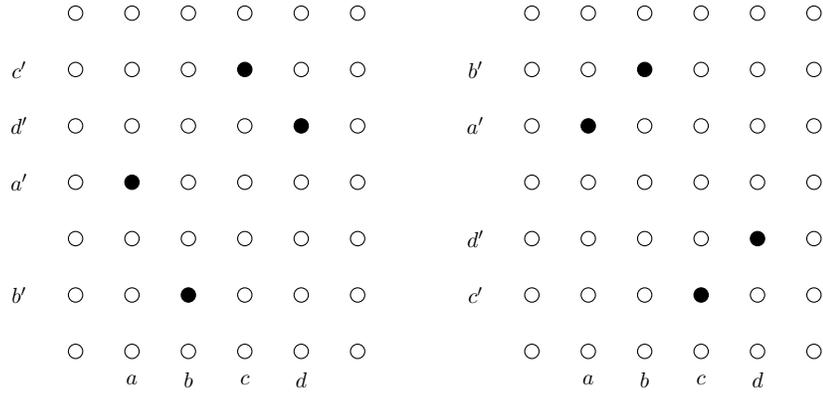

\medskip\noindent
{\bf Case 1}: $|S_1|=4$. \\
Let $S_1=\{(a',a), (b',b), (c',c), (d',d)\}$. Note that $S_1$ is then a gp-set of the grid induced by $A_1$. Since the structure of every gp-set of a grid graph is known from the proof of Theorem \ref{thm:grids-enumerate}, we can assume without loss of generality the following facts: $a < b < d$, $a < c < d$, and either ($c'< a' < b'$ and $c'< d' < b'$) or ($b'< a' < c'$ and $b'< d' < c'$). Note that there is neither relationship between $a,b$ nor between $a',b'$, namely, it can happen $a\le b$ or $b\le a$ and $a'\le b'$ or $b'\le a'$. Examples of such sets are shown in Fig.~\ref{fig-S_1}. For presentation purposes, in this and the subsequent figures an orientation is selected such that $C_s$-layers are drawn horizontally and $P_r$-layers vertically.

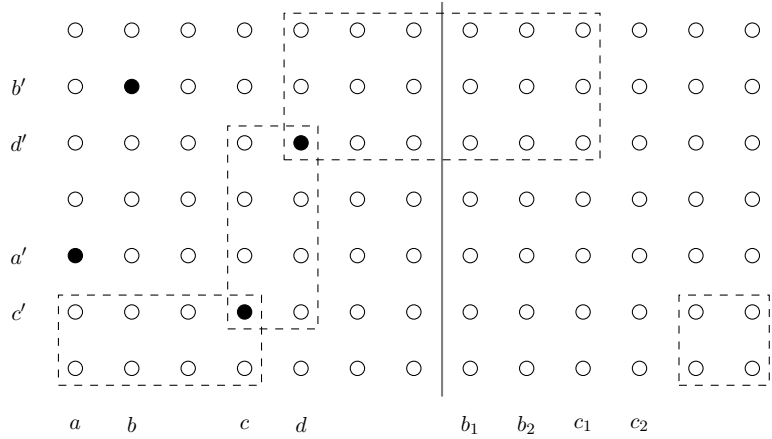
\begin{figure}[ht!]
	\centering
	\begin{tikzpicture}[scale=.75, transform shape]
	\node [draw, shape=circle,scale=0.7] (a1) at  (0,0) {};
	\node [draw, shape=circle,scale=0.7] (a2) at  (0,1) {};
	\node [draw, shape=circle,scale=0.7,fill=black] (a3) at  (0,2) {};
	\node [draw, shape=circle,scale=0.7] (a4) at  (0,3) {};
	\node [draw, shape=circle,scale=0.7] (a5) at  (0,4) {};
	\node [draw, shape=circle,scale=0.7] (a6) at  (0,5) {};
	\node [draw, shape=circle,scale=0.7] (a7) at  (0,6) {};
	
	\node [draw, shape=circle,scale=0.7] (b1) at  (1,0) {};
	\node [draw, shape=circle,scale=0.7] (b2) at  (1,1) {};
	\node [draw, shape=circle,scale=0.7] (b3) at  (1,2) {};
	\node [draw, shape=circle,scale=0.7] (b4) at  (1,3) {};
	\node [draw, shape=circle,scale=0.7] (b5) at  (1,4) {};
	\node [draw, shape=circle,scale=0.7,fill=black] (b6) at  (1,5) {};
	\node [draw, shape=circle,scale=0.7] (b7) at  (1,6) {};
	
	\node [draw, shape=circle,scale=0.7] (c1) at  (2,0) {};
	\node [draw, shape=circle,scale=0.7] (c2) at  (2,1) {};
	\node [draw, shape=circle,scale=0.7] (c3) at  (2,2) {};
	\node [draw, shape=circle,scale=0.7] (c4) at  (2,3) {};
	\node [draw, shape=circle,scale=0.7] (c5) at  (2,4) {};
	\node [draw, shape=circle,scale=0.7] (c6) at  (2,5) {};
	\node [draw, shape=circle,scale=0.7] (c7) at  (2,6) {};
	
	\node [draw, shape=circle,scale=0.7] (d1) at  (3,0) {};
	\node [draw, shape=circle,scale=0.7,fill=black] (d2) at  (3,1) {};
	\node [draw, shape=circle,scale=0.7] (d3) at  (3,2) {};
	\node [draw, shape=circle,scale=0.7] (d4) at  (3,3) {};
	\node [draw, shape=circle,scale=0.7] (d5) at  (3,4) {};
	\node [draw, shape=circle,scale=0.7] (d6) at  (3,5) {};
	\node [draw, shape=circle,scale=0.7] (d7) at  (3,6) {};
	
	\node [draw, shape=circle,scale=0.7] (e1) at  (4,0) {};
	\node [draw, shape=circle,scale=0.7] (e2) at  (4,1) {};
	\node [draw, shape=circle,scale=0.7] (e3) at  (4,2) {};
	\node [draw, shape=circle,scale=0.7] (e4) at  (4,3) {};
	\node [draw, shape=circle,scale=0.7,fill=black] (e5) at  (4,4) {};
	\node [draw, shape=circle,scale=0.7] (e6) at  (4,5) {};
	\node [draw, shape=circle,scale=0.7] (e7) at  (4,6) {};
	
	\node [draw, shape=circle,scale=0.7] (f1) at  (5,0) {};
	\node [draw, shape=circle,scale=0.7] (f2) at  (5,1) {};
	\node [draw, shape=circle,scale=0.7] (f3) at  (5,2) {};
	\node [draw, shape=circle,scale=0.7] (f4) at  (5,3) {};
	\node [draw, shape=circle,scale=0.7] (f5) at  (5,4) {};
	\node [draw, shape=circle,scale=0.7] (f6) at  (5,5) {};
	\node [draw, shape=circle,scale=0.7] (f7) at  (5,6) {};
	
	\node [draw, shape=circle,scale=0.7] (g1) at  (6,0) {};
	\node [draw, shape=circle,scale=0.7] (g2) at  (6,1) {};
	\node [draw, shape=circle,scale=0.7] (g3) at  (6,2) {};
	\node [draw, shape=circle,scale=0.7] (g4) at  (6,3) {};
	\node [draw, shape=circle,scale=0.7] (g5) at  (6,4) {};
	\node [draw, shape=circle,scale=0.7] (g6) at  (6,5) {};
	\node [draw, shape=circle,scale=0.7] (g7) at  (6,6) {};
	
	\node [draw, shape=circle,scale=0.7] (h1) at  (7,0) {};
	\node [draw, shape=circle,scale=0.7] (h2) at  (7,1) {};
	\node [draw, shape=circle,scale=0.7] (h3) at  (7,2) {};
	\node [draw, shape=circle,scale=0.7] (h4) at  (7,3) {};
	\node [draw, shape=circle,scale=0.7] (h5) at  (7,4) {};
	\node [draw, shape=circle,scale=0.7] (h6) at  (7,5) {};
	\node [draw, shape=circle,scale=0.7] (h7) at  (7,6) {};
	
	\node [draw, shape=circle,scale=0.7] (i1) at  (8,0) {};
	\node [draw, shape=circle,scale=0.7] (i2) at  (8,1) {};
	\node [draw, shape=circle,scale=0.7] (i3) at  (8,2) {};
	\node [draw, shape=circle,scale=0.7] (i4) at  (8,3) {};
	\node [draw, shape=circle,scale=0.7] (i5) at  (8,4) {};
	\node [draw, shape=circle,scale=0.7] (i6) at  (8,5) {};
	\node [draw, shape=circle,scale=0.7] (i7) at  (8,6) {};
	
	\node [draw, shape=circle,scale=0.7] (j1) at  (9,0) {};
	\node [draw, shape=circle,scale=0.7] (j2) at  (9,1) {};
	\node [draw, shape=circle,scale=0.7] (j3) at  (9,2) {};
	\node [draw, shape=circle,scale=0.7] (j4) at  (9,3) {};
	\node [draw, shape=circle,scale=0.7] (j5) at  (9,4) {};
	\node [draw, shape=circle,scale=0.7] (j6) at  (9,5) {};
	\node [draw, shape=circle,scale=0.7] (j7) at  (9,6) {};
	
	\node [draw, shape=circle,scale=0.7] (k1) at  (10,0) {};
	\node [draw, shape=circle,scale=0.7] (k2) at  (10,1) {};
	\node [draw, shape=circle,scale=0.7] (k3) at  (10,2) {};
	\node [draw, shape=circle,scale=0.7] (k4) at  (10,3) {};
	\node [draw, shape=circle,scale=0.7] (k5) at  (10,4) {};
	\node [draw, shape=circle,scale=0.7] (k6) at  (10,5) {};
	\node [draw, shape=circle,scale=0.7] (k7) at  (10,6) {};
	
	\node [draw, shape=circle,scale=0.7] (l1) at  (11,0) {};
	\node [draw, shape=circle,scale=0.7] (l2) at  (11,1) {};
	\node [draw, shape=circle,scale=0.7] (l3) at  (11,2) {};
	\node [draw, shape=circle,scale=0.7] (l4) at  (11,3) {};
	\node [draw, shape=circle,scale=0.7] (l5) at  (11,4) {};
	\node [draw, shape=circle,scale=0.7] (l6) at  (11,5) {};
	\node [draw, shape=circle,scale=0.7] (l7) at  (11,6) {};

	\node [draw, shape=circle,scale=0.7] (m1) at  (12,0) {};
	\node [draw, shape=circle,scale=0.7] (m2) at  (12,1) {};
	\node [draw, shape=circle,scale=0.7] (m3) at  (12,2) {};
	\node [draw, shape=circle,scale=0.7] (m4) at  (12,3) {};
	\node [draw, shape=circle,scale=0.7] (m5) at  (12,4) {};
	\node [draw, shape=circle,scale=0.7] (m6) at  (12,5) {};
	\node [draw, shape=circle,scale=0.7] (m7) at  (12,6) {};
	
	\draw[gray,thick] (6.5,-0.5)--(6.5,6.5);
	
	\node [scale=1] at (0,-1) {$a$};
	\node [scale=1] at (1,-1) {$b$};
	\node [scale=1] at (3,-1) {$c$};
	\node [scale=1] at (4,-1) {$d$};
	
	\node [scale=1] at (-1,1) {$c'$};
	\node [scale=1] at (-1,2) {$a'$};
	\node [scale=1] at (-1,4) {$d'$};
	\node [scale=1] at (-1,5) {$b'$};
	
	\node [scale=1] at (7,-1) {$b_1$};
	\node [scale=1] at (8,-1) {$b_2$};
	\node [scale=1] at (9,-1) {$c_1$};
	\node [scale=1] at (10,-1) {$c_2$};
	
	\draw[dashed] (3.7,3.7) rectangle (9.3,6.3);
	\draw[dashed] (2.7,0.7) rectangle (4.3,4.3);
	\draw[dashed] (-0.3,-0.3) rectangle (3.3,1.3);
	\draw[dashed] (10.7,-0.3) rectangle (12.3,1.3);
	
	\end{tikzpicture}
	\caption{The forbidden area $F((c',c),(d',d))$ appears surrounded by the dashed rectangles. Vertices $b_1$ and $b_2$ of $C_s$ are diametral with $b$ while $c_1$ and $c_2$ are diametral with $c$. Here  $b_1 = b+s_d$, $b_2 = b-s_d$, $c_1 = c+s_d$, and $c_2 = c-s_d$. Similar convention will be used in the subsequent figures. Notice that if the cycle would have an even order, then $b_1=b_2$ and $c_1=c_2$.}
	\label{forbidden-area-example}
\end{figure}

We now consider the set $F(S_1)$ in $V(P_r\cp C_s)$. Fig.~\ref{forbidden-area-example} shows an example of the forbidden area generated by only two vertices of $S_1$ ($(c',c)$ and $(d',d)$ in this case). Since it is not necessary for our purposes, we do not look at the whole set $S_1$, only a significant part.

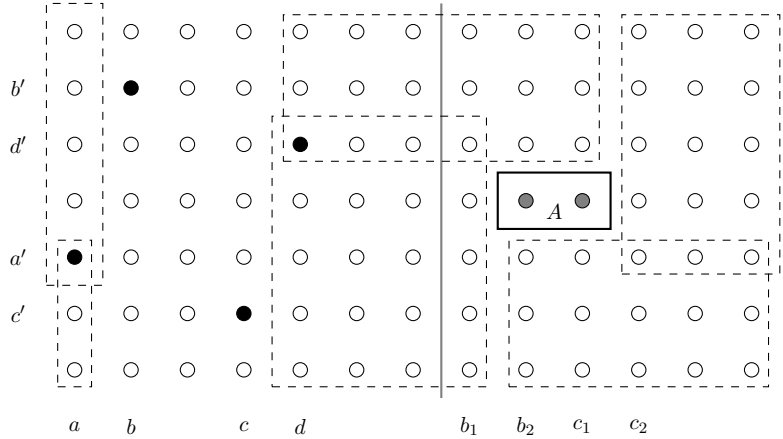
\begin{figure}[ht!]
	\centering
	\begin{tikzpicture}[scale=.75, transform shape]
	\node [draw, shape=circle,scale=0.7] (a1) at  (0,0) {};
	\node [draw, shape=circle,scale=0.7] (a2) at  (0,1) {};
	\node [draw, shape=circle,scale=0.7,fill=black] (a3) at  (0,2) {};
	\node [draw, shape=circle,scale=0.7] (a4) at  (0,3) {};
	\node [draw, shape=circle,scale=0.7] (a5) at  (0,4) {};
	\node [draw, shape=circle,scale=0.7] (a6) at  (0,5) {};
	\node [draw, shape=circle,scale=0.7] (a7) at  (0,6) {};
	
	\node [draw, shape=circle,scale=0.7] (b1) at  (1,0) {};
	\node [draw, shape=circle,scale=0.7] (b2) at  (1,1) {};
	\node [draw, shape=circle,scale=0.7] (b3) at  (1,2) {};
	\node [draw, shape=circle,scale=0.7] (b4) at  (1,3) {};
	\node [draw, shape=circle,scale=0.7] (b5) at  (1,4) {};
	\node [draw, shape=circle,scale=0.7,fill=black] (b6) at  (1,5) {};
	\node [draw, shape=circle,scale=0.7] (b7) at  (1,6) {};
	
	\node [draw, shape=circle,scale=0.7] (c1) at  (2,0) {};
	\node [draw, shape=circle,scale=0.7] (c2) at  (2,1) {};
	\node [draw, shape=circle,scale=0.7] (c3) at  (2,2) {};
	\node [draw, shape=circle,scale=0.7] (c4) at  (2,3) {};
	\node [draw, shape=circle,scale=0.7] (c5) at  (2,4) {};
	\node [draw, shape=circle,scale=0.7] (c6) at  (2,5) {};
	\node [draw, shape=circle,scale=0.7] (c7) at  (2,6) {};
	
	\node [draw, shape=circle,scale=0.7] (d1) at  (3,0) {};
	\node [draw, shape=circle,scale=0.7,fill=black] (d2) at  (3,1) {};
	\node [draw, shape=circle,scale=0.7] (d3) at  (3,2) {};
	\node [draw, shape=circle,scale=0.7] (d4) at  (3,3) {};
	\node [draw, shape=circle,scale=0.7] (d5) at  (3,4) {};
	\node [draw, shape=circle,scale=0.7] (d6) at  (3,5) {};
	\node [draw, shape=circle,scale=0.7] (d7) at  (3,6) {};
	
	\node [draw, shape=circle,scale=0.7] (e1) at  (4,0) {};
	\node [draw, shape=circle,scale=0.7] (e2) at  (4,1) {};
	\node [draw, shape=circle,scale=0.7] (e3) at  (4,2) {};
	\node [draw, shape=circle,scale=0.7] (e4) at  (4,3) {};
	\node [draw, shape=circle,scale=0.7,fill=black] (e5) at  (4,4) {};
	\node [draw, shape=circle,scale=0.7] (e6) at  (4,5) {};
	\node [draw, shape=circle,scale=0.7] (e7) at  (4,6) {};
	
	\node [draw, shape=circle,scale=0.7] (f1) at  (5,0) {};
	\node [draw, shape=circle,scale=0.7] (f2) at  (5,1) {};
	\node [draw, shape=circle,scale=0.7] (f3) at  (5,2) {};
	\node [draw, shape=circle,scale=0.7] (f4) at  (5,3) {};
	\node [draw, shape=circle,scale=0.7] (f5) at  (5,4) {};
	\node [draw, shape=circle,scale=0.7] (f6) at  (5,5) {};
	\node [draw, shape=circle,scale=0.7] (f7) at  (5,6) {};
	
	\node [draw, shape=circle,scale=0.7] (g1) at  (6,0) {};
	\node [draw, shape=circle,scale=0.7] (g2) at  (6,1) {};
	\node [draw, shape=circle,scale=0.7] (g3) at  (6,2) {};
	\node [draw, shape=circle,scale=0.7] (g4) at  (6,3) {};
	\node [draw, shape=circle,scale=0.7] (g5) at  (6,4) {};
	\node [draw, shape=circle,scale=0.7] (g6) at  (6,5) {};
	\node [draw, shape=circle,scale=0.7] (g7) at  (6,6) {};
	
	\node [draw, shape=circle,scale=0.7] (h1) at  (7,0) {};
	\node [draw, shape=circle,scale=0.7] (h2) at  (7,1) {};
	\node [draw, shape=circle,scale=0.7] (h3) at  (7,2) {};
	\node [draw, shape=circle,scale=0.7] (h4) at  (7,3) {};
	\node [draw, shape=circle,scale=0.7] (h5) at  (7,4) {};
	\node [draw, shape=circle,scale=0.7] (h6) at  (7,5) {};
	\node [draw, shape=circle,scale=0.7] (h7) at  (7,6) {};
	
	\node [draw, shape=circle,scale=0.7] (i1) at  (8,0) {};
	\node [draw, shape=circle,scale=0.7] (i2) at  (8,1) {};
	\node [draw, shape=circle,scale=0.7] (i3) at  (8,2) {};
	\node [draw, shape=circle,scale=0.7,fill=gray] (i4) at  (8,3) {};
	\node [draw, shape=circle,scale=0.7] (i5) at  (8,4) {};
	\node [draw, shape=circle,scale=0.7] (i6) at  (8,5) {};
	\node [draw, shape=circle,scale=0.7] (i7) at  (8,6) {};
	
	\node [draw, shape=circle,scale=0.7] (j1) at  (9,0) {};
	\node [draw, shape=circle,scale=0.7] (j2) at  (9,1) {};
	\node [draw, shape=circle,scale=0.7] (j3) at  (9,2) {};
	\node [draw, shape=circle,scale=0.7,fill=gray] (j4) at  (9,3) {};
	\node [draw, shape=circle,scale=0.7] (j5) at  (9,4) {};
	\node [draw, shape=circle,scale=0.7] (j6) at  (9,5) {};
	\node [draw, shape=circle,scale=0.7] (j7) at  (9,6) {};
	
	\node [draw, shape=circle,scale=0.7] (k1) at  (10,0) {};
	\node [draw, shape=circle,scale=0.7] (k2) at  (10,1) {};
	\node [draw, shape=circle,scale=0.7] (k3) at  (10,2) {};
	\node [draw, shape=circle,scale=0.7] (k4) at  (10,3) {};
	\node [draw, shape=circle,scale=0.7] (k5) at  (10,4) {};
	\node [draw, shape=circle,scale=0.7] (k6) at  (10,5) {};
	\node [draw, shape=circle,scale=0.7] (k7) at  (10,6) {};
	
	\node [draw, shape=circle,scale=0.7] (l1) at  (11,0) {};
	\node [draw, shape=circle,scale=0.7] (l2) at  (11,1) {};
	\node [draw, shape=circle,scale=0.7] (l3) at  (11,2) {};
	\node [draw, shape=circle,scale=0.7] (l4) at  (11,3) {};
	\node [draw, shape=circle,scale=0.7] (l5) at  (11,4) {};
	\node [draw, shape=circle,scale=0.7] (l6) at  (11,5) {};
	\node [draw, shape=circle,scale=0.7] (l7) at  (11,6) {};

	\node [draw, shape=circle,scale=0.7] (m1) at  (12,0) {};
	\node [draw, shape=circle,scale=0.7] (m2) at  (12,1) {};
	\node [draw, shape=circle,scale=0.7] (m3) at  (12,2) {};
	\node [draw, shape=circle,scale=0.7] (m4) at  (12,3) {};
	\node [draw, shape=circle,scale=0.7] (m5) at  (12,4) {};
	\node [draw, shape=circle,scale=0.7] (m6) at  (12,5) {};
	\node [draw, shape=circle,scale=0.7] (m7) at  (12,6) {};
	
	\draw[gray,thick] (6.5,-0.5)--(6.5,6.5);
	
	\node [scale=1] at (0,-1) {$a$};
	\node [scale=1] at (1,-1) {$b$};
	\node [scale=1] at (3,-1) {$c$};
	\node [scale=1] at (4,-1) {$d$};
	
	\node [scale=1] at (-1,1) {$c'$};
	\node [scale=1] at (-1,2) {$a'$};
	\node [scale=1] at (-1,4) {$d'$};
	\node [scale=1] at (-1,5) {$b'$};
	
	\node [scale=1] at (7,-1) {$b_1$};
	\node [scale=1] at (8,-1) {$b_2$};
	\node [scale=1] at (9,-1) {$c_1$};
	\node [scale=1] at (10,-1) {$c_2$};
	
	\draw[dashed] (3.7,3.7) rectangle (9.3,6.3);
	\draw[dashed] (3.5,4.5) rectangle (7.3,-0.3);
	\draw[dashed] (-0.3,-0.3) rectangle (0.3,2.3);
	\draw[dashed] (-0.5,1.5) rectangle (0.5,6.5);
	\draw[dashed] (7.7,-0.3) rectangle (12.3,2.3);
	\draw[dashed] (12.5,1.7) rectangle (9.7,6.3);
	
	\draw[thick] (7.5,2.5) rectangle (9.5,3.5);
	
	\node [scale=1] at (8.5,2.8) {$A$};
	
	\end{tikzpicture}
	\caption{Part of the forbidden area of the bolded set of vertices appears in dashed rectangles. The two gray vertices of the thick rectangle (denoted by $A$) do not belong to the forbidden area of the bolded vertices.}
	\label{forbidden-area-set}
\end{figure}

We detail now the case $a < b \le c < d$ and $c'< a' \le d' < b'$, see Fig.~\ref{forbidden-area-set}. Observe that:
\begin{align*}
\{a',\dots,r-1\}\times \{a,a-1,\dots,c-s_d\}& \subset F((c',c),(a',a)), \\
\{0,\dots,a'\}\times \{a,a-1,\dots,b-s_d\}& \subset F((b',b),(a',a)), \\
\{d',\dots,r-1\}\times \{d,d+1,\dots,c+s_d\}& \subset F((c',c),(d',d)), \\
\{0,\dots,d'\}\times \{d,d+1,\dots,b+s_d\}& \subset F((b',b),(d',d)).
\end{align*}
Let us define a subset $A$ of $A_2$ by $A=A_2-F(S_1)$. Observe that $A$ is empty if and only if one of the following situations occur: $a'=d'$, or $a'=d'-1$, or $b=c$, or ($b=c-1$ and $s$ is even). More precisely, $A=\{a'+1,\dots,d'-1\}\times \{b_2,\dots,c_1\}$. Fig.~\ref{forbidden-area-set} shows an example of this where the set $A$ is not empty.

As a consequence, it must happen that $S_2\subseteq A$, since otherwise we get a contradiction with $S$ being a gp-set. If $|S_2|\ge 2$, then let $(x',x),(y',y)\in S_2$. It is then not difficult to observe that $(c',c),(x',x),(y',y)$ or $(d',d),(x',x),(y',y)$ lie in a geodesic of $P_r\cp C_s$, as well as $(a',a),(x',x),(y',y)$ or $(b',b),(x',x),(y',y)$, which is not possible. Thus $|S_2|\le 1$.

By using a similar reasoning, we deduce the same conclusion for any other relationship of $a,b,c,d$ and $a',b',c',d'$ (from those ones that allow to obtain a gp-set of the grid induced by $A_1$, according to the proof of Theorem \ref{thm:grids-enumerate}). As a consequence of the whole deduction of this case, we obtain that $\gp(P_r\cp C_s)=|S|=|S_1|+|S_2|\le 5$.

\medskip\noindent
{\bf Case 2}: $|S_2|=4$. \\
The situation is similar to Case 1, although if $s$ is odd, then the set $A_2$ is smaller than $A_1$ by a difference of one $P_r$-layer. Notwithstanding, this fact does not influence the arguments considered in Case 1.

\medskip\noindent
{\bf Case 3}: $|S_1|=3$ or $|S_2|=3$.\\
Assume first that $|S_1|=3$, and let $S_1=\{(a',a),(b',b),(c',c)\}$. Clearly, the three elements of $S_1$ cannot lie simultaneously in the same $^iC_s$-layer, or in the same $P_r\,^j$-layer. Moreover, it cannot happen that $a'\le b'\le c'$ and $a\le b\le c$ at the same time, or any other similar double monotone sequence. This means that, for instance, if $a'\le b'\le c'$, then either ($b < a$ and $b < c$) or ($b > a$ and $b > c$).

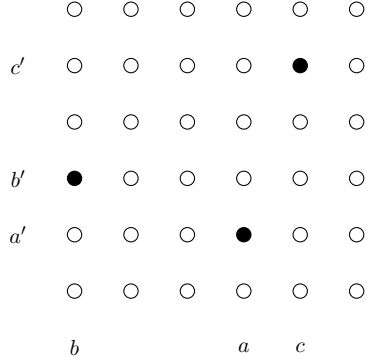
\begin{figure}[ht!]
	\centering
	\begin{tikzpicture}[scale=.75, transform shape]
	\node [draw, shape=circle,scale=0.7] (a1) at  (0,0) {};
	\node [draw, shape=circle,scale=0.7] (a2) at  (0,1) {};
	\node [draw, shape=circle,scale=0.7,fill=black] (a3) at  (0,2) {};
	\node [draw, shape=circle,scale=0.7] (a4) at  (0,3) {};
	\node [draw, shape=circle,scale=0.7] (a5) at  (0,4) {};
	\node [draw, shape=circle,scale=0.7] (a6) at  (0,5) {};
	
	\node [draw, shape=circle,scale=0.7] (b1) at  (1,0) {};
	\node [draw, shape=circle,scale=0.7] (b2) at  (1,1) {};
	\node [draw, shape=circle,scale=0.7] (b3) at  (1,2) {};
	\node [draw, shape=circle,scale=0.7] (b4) at  (1,3) {};
	\node [draw, shape=circle,scale=0.7] (b5) at  (1,4) {};
	\node [draw, shape=circle,scale=0.7] (b6) at  (1,5) {};
	
	\node [draw, shape=circle,scale=0.7] (c1) at  (2,0) {};
	\node [draw, shape=circle,scale=0.7] (c2) at  (2,1) {};
	\node [draw, shape=circle,scale=0.7] (c3) at  (2,2) {};
	\node [draw, shape=circle,scale=0.7] (c4) at  (2,3) {};
	\node [draw, shape=circle,scale=0.7] (c5) at  (2,4) {};
	\node [draw, shape=circle,scale=0.7] (c6) at  (2,5) {};
	
	\node [draw, shape=circle,scale=0.7] (d1) at  (3,0) {};
	\node [draw, shape=circle,scale=0.7,fill=black] (d2) at  (3,1) {};
	\node [draw, shape=circle,scale=0.7] (d3) at  (3,2) {};
	\node [draw, shape=circle,scale=0.7] (d4) at  (3,3) {};
	\node [draw, shape=circle,scale=0.7] (d5) at  (3,4) {};
	\node [draw, shape=circle,scale=0.7] (d6) at  (3,5) {};
	
	\node [draw, shape=circle,scale=0.7] (e1) at  (4,0) {};
	\node [draw, shape=circle,scale=0.7] (e2) at  (4,1) {};
	\node [draw, shape=circle,scale=0.7] (e3) at  (4,2) {};
	\node [draw, shape=circle,scale=0.7] (e4) at  (4,3) {};
	\node [draw, shape=circle,scale=0.7,fill=black] (e5) at  (4,4) {};
	\node [draw, shape=circle,scale=0.7] (e6) at  (4,5) {};
	
	\node [draw, shape=circle,scale=0.7] (f1) at  (5,0) {};
	\node [draw, shape=circle,scale=0.7] (f2) at  (5,1) {};
	\node [draw, shape=circle,scale=0.7] (f3) at  (5,2) {};
	\node [draw, shape=circle,scale=0.7] (f4) at  (5,3) {};
	\node [draw, shape=circle,scale=0.7] (f5) at  (5,4) {};
	\node [draw, shape=circle,scale=0.7] (f6) at  (5,5) {};
	
	\node [scale=1] at (0,-1) {$b$};
	\node [scale=1] at (3,-1) {$a$};
	\node [scale=1] at (4,-1) {$c$};
	
	\node [scale=1] at (-1,1) {$a'$};
	\node [scale=1] at (-1,2) {$b'$};
	\node [scale=1] at (-1,4) {$c'$};

	\end{tikzpicture}
	\caption{An example of a configuration for the set $S_1$.}
	\label{fig-S_1-case-2}
\end{figure}

We may assume now that $a'\le b'\le c'$, $b < a$ and $b < c$. Fig.~\ref{fig-S_1-case-2} shows an example of this.

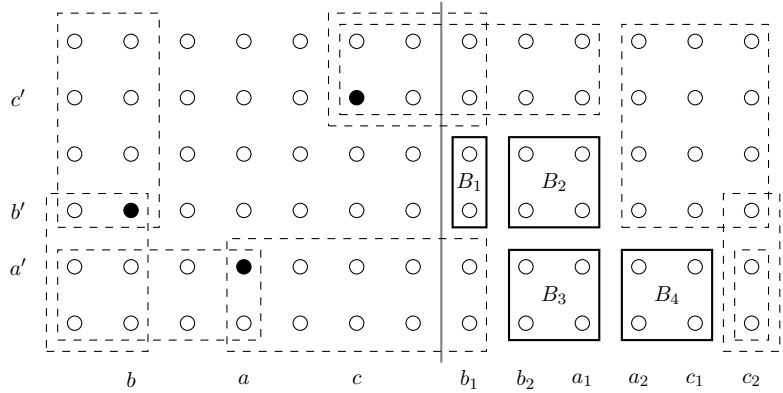
\begin{figure}[ht!]
	\centering
	\begin{tikzpicture}[scale=.75, transform shape]
	\node [draw, shape=circle,scale=0.7] (a1) at  (0,0) {};
	\node [draw, shape=circle,scale=0.7] (a2) at  (0,1) {};
	\node [draw, shape=circle,scale=0.7] (a3) at  (0,2) {};
	\node [draw, shape=circle,scale=0.7] (a4) at  (0,3) {};
	\node [draw, shape=circle,scale=0.7] (a5) at  (0,4) {};
	\node [draw, shape=circle,scale=0.7] (a6) at  (0,5) {};
	
	\node [draw, shape=circle,scale=0.7] (b1) at  (1,0) {};
	\node [draw, shape=circle,scale=0.7] (b2) at  (1,1) {};
	\node [draw, shape=circle,scale=0.7,fill=black] (b3) at  (1,2) {};
	\node [draw, shape=circle,scale=0.7] (b4) at  (1,3) {};
	\node [draw, shape=circle,scale=0.7] (b5) at  (1,4) {};
	\node [draw, shape=circle,scale=0.7] (b6) at  (1,5) {};
	
	\node [draw, shape=circle,scale=0.7] (c1) at  (2,0) {};
	\node [draw, shape=circle,scale=0.7] (c2) at  (2,1) {};
	\node [draw, shape=circle,scale=0.7] (c3) at  (2,2) {};
	\node [draw, shape=circle,scale=0.7] (c4) at  (2,3) {};
	\node [draw, shape=circle,scale=0.7] (c5) at  (2,4) {};
	\node [draw, shape=circle,scale=0.7] (c6) at  (2,5) {};
	
	\node [draw, shape=circle,scale=0.7] (d1) at  (3,0) {};
	\node [draw, shape=circle,scale=0.7,fill=black] (d2) at  (3,1) {};
	\node [draw, shape=circle,scale=0.7] (d3) at  (3,2) {};
	\node [draw, shape=circle,scale=0.7] (d4) at  (3,3) {};
	\node [draw, shape=circle,scale=0.7] (d5) at  (3,4) {};
	\node [draw, shape=circle,scale=0.7] (d6) at  (3,5) {};
	
	\node [draw, shape=circle,scale=0.7] (e1) at  (4,0) {};
	\node [draw, shape=circle,scale=0.7] (e2) at  (4,1) {};
	\node [draw, shape=circle,scale=0.7] (e3) at  (4,2) {};
	\node [draw, shape=circle,scale=0.7] (e4) at  (4,3) {};
	\node [draw, shape=circle,scale=0.7] (e5) at  (4,4) {};
	\node [draw, shape=circle,scale=0.7] (e6) at  (4,5) {};
	
	\node [draw, shape=circle,scale=0.7] (f1) at  (5,0) {};
	\node [draw, shape=circle,scale=0.7] (f2) at  (5,1) {};
	\node [draw, shape=circle,scale=0.7] (f3) at  (5,2) {};
	\node [draw, shape=circle,scale=0.7] (f4) at  (5,3) {};
	\node [draw, shape=circle,scale=0.7,fill=black] (f5) at  (5,4) {};
	\node [draw, shape=circle,scale=0.7] (f6) at  (5,5) {};
	
	\node [draw, shape=circle,scale=0.7] (g1) at  (6,0) {};
	\node [draw, shape=circle,scale=0.7] (g2) at  (6,1) {};
	\node [draw, shape=circle,scale=0.7] (g3) at  (6,2) {};
	\node [draw, shape=circle,scale=0.7] (g4) at  (6,3) {};
	\node [draw, shape=circle,scale=0.7] (g5) at  (6,4) {};
	\node [draw, shape=circle,scale=0.7] (g6) at  (6,5) {};
	
	\node [draw, shape=circle,scale=0.7] (h1) at  (7,0) {};
	\node [draw, shape=circle,scale=0.7] (h2) at  (7,1) {};
	\node [draw, shape=circle,scale=0.7] (h3) at  (7,2) {};
	\node [draw, shape=circle,scale=0.7] (h4) at  (7,3) {};
	\node [draw, shape=circle,scale=0.7] (h5) at  (7,4) {};
	\node [draw, shape=circle,scale=0.7] (h6) at  (7,5) {};
	
	\node [draw, shape=circle,scale=0.7] (i1) at  (8,0) {};
	\node [draw, shape=circle,scale=0.7] (i2) at  (8,1) {};
	\node [draw, shape=circle,scale=0.7] (i3) at  (8,2) {};
	\node [draw, shape=circle,scale=0.7] (i4) at  (8,3) {};
	\node [draw, shape=circle,scale=0.7] (i5) at  (8,4) {};
	\node [draw, shape=circle,scale=0.7] (i6) at  (8,5) {};
	
	\node [draw, shape=circle,scale=0.7] (j1) at  (9,0) {};
	\node [draw, shape=circle,scale=0.7] (j2) at  (9,1) {};
	\node [draw, shape=circle,scale=0.7] (j3) at  (9,2) {};
	\node [draw, shape=circle,scale=0.7] (j4) at  (9,3) {};
	\node [draw, shape=circle,scale=0.7] (j5) at  (9,4) {};
	\node [draw, shape=circle,scale=0.7] (j6) at  (9,5) {};
	
	\node [draw, shape=circle,scale=0.7] (k1) at  (10,0) {};
	\node [draw, shape=circle,scale=0.7] (k2) at  (10,1) {};
	\node [draw, shape=circle,scale=0.7] (k3) at  (10,2) {};
	\node [draw, shape=circle,scale=0.7] (k4) at  (10,3) {};
	\node [draw, shape=circle,scale=0.7] (k5) at  (10,4) {};
	\node [draw, shape=circle,scale=0.7] (k6) at  (10,5) {};
	
	\node [draw, shape=circle,scale=0.7] (l1) at  (11,0) {};
	\node [draw, shape=circle,scale=0.7] (l2) at  (11,1) {};
	\node [draw, shape=circle,scale=0.7] (l3) at  (11,2) {};
	\node [draw, shape=circle,scale=0.7] (l4) at  (11,3) {};
	\node [draw, shape=circle,scale=0.7] (l5) at  (11,4) {};
	\node [draw, shape=circle,scale=0.7] (l6) at  (11,5) {};
	
	\node [draw, shape=circle,scale=0.7] (m1) at  (12,0) {};
	\node [draw, shape=circle,scale=0.7] (m2) at  (12,1) {};
	\node [draw, shape=circle,scale=0.7] (m3) at  (12,2) {};
	\node [draw, shape=circle,scale=0.7] (m4) at  (12,3) {};
	\node [draw, shape=circle,scale=0.7] (m5) at  (12,4) {};
	\node [draw, shape=circle,scale=0.7] (m6) at  (12,5) {};
	
	\node [scale=1] at (1,-1) {$b$};
	\node [scale=1] at (3,-1) {$a$};
	\node [scale=1] at (5,-1) {$c$};
	
	\node [scale=1] at (-1,1) {$a'$};
	\node [scale=1] at (-1,2) {$b'$};
	\node [scale=1] at (-1,4) {$c'$};
	
	\node [scale=1] at (7,-1) {$b_1$};
	\node [scale=1] at (8,-1) {$b_2$};
	\node [scale=1] at (9,-1) {$a_1$};
	\node [scale=1] at (10,-1) {$a_2$};
	\node [scale=1] at (11,-1) {$c_1$};
	\node [scale=1] at (12,-1) {$c_2$};
	
	\draw[gray,thick] (6.5,-0.7)--(6.5,5.7);
	
	\draw[dashed] (4.7,3.7) rectangle (9.3,5.3);
	\draw[dashed] (-0.3,-0.3) rectangle (3.3,1.3);
	\draw[dashed] (11.7,-0.3) rectangle (12.3,1.3);
	\draw[dashed] (4.5,3.5) rectangle (7.3,5.5);
	\draw[dashed] (-0.5,-0.5) rectangle (1.3,2.3);
	\draw[dashed] (11.5,-0.5) rectangle (12.5,2.3);
	\draw[dashed] (2.7,-0.5) rectangle (7.3,1.5);
	\draw[dashed] (-0.3,1.7) rectangle (1.5,5.5);
	\draw[dashed] (9.7,1.7) rectangle (12.3,5.3);
	
	\draw[thick] (6.7,1.7) rectangle (7.3,3.3);
	\draw[thick] (7.7,-0.3) rectangle (9.3,1.3);
	\draw[thick] (7.7,1.7) rectangle (9.3,3.3);
	\draw[thick] (9.7,-0.3) rectangle (11.3,1.3);
	
	\node [scale=1] at (7,2.5) {$B_1$};
	\node [scale=1] at (8.5,2.5) {$B_2$};
	\node [scale=1] at (8.5,0.5) {$B_3$};
	\node [scale=1] at (10.5,0.5) {$B_4$};
	
	\end{tikzpicture}
	\caption{A significant part of the set $F(S_1)$ appears surrounded by dashed rectangles. For $x\in \{a,b,c\}$, the vertices $x_1$ and $x_2$ from $C_s$ ($x_1=x+s_d$ and $x_2=x-s_d$) are diametral vertices with $x$. Note that if $C_s$ is an even cycle, then $x_1 = x_2$.}
	\label{fig-S_1-case-2-F}
\end{figure}

We now consider the set $F(S_1)$ in $V(P_r\cp C_s)$, and observe the following. Recalling that $s_d=\left\lfloor s/2\right\rfloor$ we have:
\begin{align*}
\{c',\dots,r-1\}\times \{c,c+1,\dots,a+s_d\} & \subset F((c',c),(a',a)), \\
\{0,\dots,a'\}\times\{a,a-1,\dots,c-s_d\} & \subset F((c',c),(a',a)),\\
\{b',\dots,r-1\}\times \{b,b-1,\dots,a-s_d\} & \subset F((b',b),(a',a)), \\
\{0,\dots,a'\}\times \{a,a+1,\dots,b+s_d\} & \subset F((b',b),(a',a)), \\
\{c',\dots,r-1\}\times \{c,c+1,\dots,b+s_d\} & \subset F((c',c),(b',b)), \\
\{0,\dots,b'\}\times\{b,b-1,\dots,c-s_d\} & \subset F((c',c),(b',b)).
\end{align*}
See Fig.~\ref{fig-S_1-case-2-F} for an example of the situations above.

Observe now that there are four sets, say $B_1$, $B_2$, $B_3$, and $B_4$, such that $B_1\cup B_2\cup B_3\cup B_4=A_2-F(S_1)$, and satisfying the following. If $B_1$, $B_2$, $B_3$, and $B_4$ are not empty, then
\begin{align*}
B_1 & = \{a'+1,\dots,c'-1\}\times\{s_d+1,\dots,b+s_d\}, \\
B_2 & = \{a'+1,\dots,c'-1\}\times\{b+s_d+1,\dots,a+s_d\}, \\
B_3 & = \{0,\dots,a'\}\times\{b+s_d+1,\dots,a+s_d\}, \\
B_4 & = \{0,\dots,a'\}\times\{a+s_d+1,\dots,c+s_d\}.
\end{align*}
Note that some of these sets could be empty, or could have non-empty intersection, depending on the parity of $s$ and on the structure of the set $S_1$.

If $|S_2\cap B_i|\ge 2$ for some $i\in [4]$, then we shall find an isometric subgraph of $P_r\cp C_s$ isomorphic to a grid graph such that it contains four vertices of the set $S$. Hence, we can change the partition given by $A_1$ and $A_2$ from the beginning, to a new one, and proceed as in Case 1, to prove that $\gp(P_r\cp C_s)\le 5$. That is, if $|S_2\cap B_i|\ge 2$ for some $i\in [3]$, then we can use the partition $A'_1=[r]_0\times \{a,a+1,\dots,a+s_d\}$ and $A'_2=V(P_r\cp C_s)-A'_1$, and if $|S_2\cap B_4|\ge 2$, then we can use the partition $A'_1=[r]_0\times \{a+s_d+1,a+s_d+2,\dots,a\}$ (note that $a+s_d+1=a-s_d$) and $A'_2=V(P_r\cp C_s)-A'_1$. In concordance, we may assume that $|S_2\cap B_i|\le 1$ for every $i\in [4]$.

We consider now the three sets $B_1$, $B_2$ and $B_3$. If at least two of them contain one element from $S_2$, then, as above, we can find a different partition of $V(P_r\cp C_s)$ and proceed like in Case 1. Thus, $|(B_1\cup B_2\cup B_3)\cap S_1|\le 1$.

Finally, we deduce that $\gp(P_r\cp C_s)=|S|=|S_1|+|S_2|= |S_1|+ |(B_1\cup B_2\cup B_3)\cap S_1|+|B_4\cap S_1|\le 5$. By using similar arguments, we can again obtain a similar conclusion for any other possible relationship between $a,b,c$ and $a',b',c'$.

Now, if $|S_2|=3$ and $s$ is even, then an identical procedure as above leads to the same conclusion. Finally, if $|S_2|=3$ and $s$ is odd, then the arguments above do not apply since $A_2$ has only $\lfloor s/2\rfloor$ $P_r$-layers. However, it must happen $|S_1|\le 2$, since otherwise, we proceed as in the previous situation (when $|S_1|=3$). Consequently, we have that $\gp(P_r\cp C_s)=|S|=|S_1|+|S_2|\le 5$. This completes the proof of this case, and therefore, of the whole theorem.
\qed

\section{Torus graphs}
\label{sec:torus}

Knowing $\gp(P_r\cp P_s)$ and $\gp(P_r\cp C_s)$, the next task is to consider the {\em torus graphs} $C_r\cp C_s$, $r, s\ge 3$, where we keep the convention that $V(C_n) = [n]_0$. In contrast to the former two cases, for the torus graphs we are not able to give an exact result, but we will prove a theorem which reduces the general position number of the torus graphs to only two cases. The following lemma and short remarks will be useful for our proof.

\begin{lemma}
\label{lem:grid}
If $S$ is a general position set in $P_r \cp P_s$ and there exists $x \in S$, with ${\rm deg}(x)=2$, then $|S|\leq 3.$
\end{lemma}
\proof
Suppose on the contrary that $|S|=4$. Assume without loss of generality that $(0,0)\in S.$ Let $(i,i'),(j,j'),(k,k')$ be the other three vertices of $S$ where $i\leq j \leq k$. If $i'\leq j'$, then $(0,0),(i,i'),(j,j')$ lie on a common geodesic. So $i'>j'.$ Similarly, $j'>k'.$ But now $(i,i'),(j,j'),(k,k')$ lie on a common geodesic.
\qed
\begin{remark}\label{rem:different}
If $S = \{(0,0),(i,i'),(j,j')\}$ is a general position set of $P_r \cp P_s$, then $i\neq j$ and $i' \neq j'$.
\end{remark}
\begin{remark}\label{rem:isometric}
The graph $P_r$ is an isometric subgraph of $C_s$ if $r \leq \lfloor \frac{s}{2}\rfloor +1.$
\end{remark}

We now bound $\gp(C_r\cp C_s)$ from the above as follows. 

\begin{theorem}
\label{thm:torus_up}
If $r\ge 3$ and $s\ge 3$, then $\gp(C_r\cp C_s)\leq 7$.
\end{theorem}

\proof
Let $S\subseteq V(C_r \cp C_s)$ be a general position set and let $w\in S$. Since $C_r \cp C_s$ is a vertex-transitive graph (meaning that for each pair of vertices of $C_r \cp C_s$ there is an automorphism which maps the first vertex into the second), we can without loss of generality set $w =  (\lfloor \frac{r}{2} \rfloor, \lfloor \frac{s}{2}\rfloor)$. Consider the following four subgraphs of  $C_r \cp C_s$ (where addition is done modulo $r$ or $s$):
\begin{itemize}
    \item $X_1 = \{(u_i,v_j): u_i \in [\lfloor \frac{r}{2}\rfloor+1]_0, v_j \in [\lfloor \frac{s}{2}\rfloor+1]_0\},$
    \item $X_2 = \{(u_i,v_j): u_i \in \{\lfloor \frac{r}{2}\rfloor,\lfloor \frac{r}{2}\rfloor +1,\ldots,\lfloor \frac{r}{2}\rfloor+\lfloor \frac{r}{2}\rfloor\}, v_j \in [\lfloor \frac{s}{2}\rfloor+1]_0\},$
    \item $X_3 = \{(u_i,v_j): u_i \in \{\lfloor \frac{r}{2}\rfloor,\lfloor \frac{r}{2}\rfloor +1,\ldots,\lfloor \frac{r}{2}\rfloor+\lfloor \frac{r}{2}\rfloor\}, v_j \in \{\lfloor \frac{s}{2}\rfloor,\lfloor \frac{s}{2}\rfloor +1,\ldots,\lfloor \frac{s}{2}\rfloor+\lfloor \frac{s}{2}\rfloor\},$
    \item $X_4 = \{(u_i,v_j): u_i \in [\lfloor \frac{r}{2}\rfloor+1]_0, v_j \in \{\lfloor \frac{s}{2}\rfloor,\lfloor \frac{s}{2}\rfloor +1,\ldots,\lfloor \frac{s}{2}\rfloor+\lfloor \frac{s}{2}\rfloor\}.$
\end{itemize}

For each $i \in [4]$, the graph $X_i$ is isomorphic to $P_{\lfloor\frac{r}{2}\rfloor +1} \cp P_{\lfloor\frac{s}{2}\rfloor +1}$. Hence, by Remark \ref{rem:isometric}, every $X_i$ is an isometric subgraph of $C_r \cp C_s.$ Since $w \in X_i$, $i \in [4]$, it is clear that $\gp(C_r \cp C_s) \leq 9$ (by Lemma \ref{lem:grid}). This bound can be further improved using the following:

\medskip\noindent
\textbf{Claim}: If $|(X_1 \cup X_3) \cap S| = 5$, then $|(X_2 \cup X_4) \cap S| \leq 2$.

\medskip\noindent
Clearly, if the claim is proved, then by symmetry we also infer that if $|(X_2 \cup X_4) \cap S| = 5$, then $|(X_1 \cup X_3) \cap S| \leq 2$. From these two facts the assertion of the theorem will follow. It thus remains to prove the claim. 

Let $u_1 = (x_1,y_1)$, $u_2 = (x_2,y_2)$, $u_3 = (x_3,y_3)$, and $u_4 = (x_4,y_4)$ be vertices such that $X_1 \cap S =  \{u_1,u_2,w\}$ and $X_3 \cap S =\{u_3,u_4,w\}$. Suppose that there are vertices $u_5 = (x_5,y_5)$ and $u_6 = (x_6,y_6)$ such that $\{u_5,u_6\} \subseteq X_2 \cap S$. As $S$ is a general position, Remark~\ref{rem:different} implies that we may without loss of generality assume that $x_1<x_2<x_3<x_4$ and $y_2<y_1<y_4<y_3$. Further, we may also assume that $x_5<x_6$ and $y_5<y_6$.

If $x_3 - x_1 \leq \frac{r}{2}$, then either $w$ lies on a $u_1,u_3$-geodesic  (if $y_3-y_1 \leq \frac{s}{2}$) or $u_2$ lies on such a geodesic (when $y_3 - y_1 \geq \frac{s}{2}$). It follows that $x_3 - x_1 > \frac{r}{2}$. With a  similar reasoning we also get that $y_3 - y_1 > \frac{s}{2}$, since otherwise $u_4$ is on a $u_1,u_3$-geodesic.

Applying the same reasoning the following is obtained: $x_4 - x_2 > \frac{r}{2}$ and $y_4 - y_2 > \frac{s}{2}$. We now distinguish the following cases. 

\medskip\noindent
{\bf Case 1}: ($y_3 - y_6 \geq \frac{s}{2}$).\\
Subcase 1a: $x_3 < x_6$.\\
First, $x_3>x_5$, since otherwise $u_5$ lies on a $u_3,u_6$-geodesic. If $y_4-y_5 \geq \frac{s}{2}$, then $u_3$ lies on a $u_4,u_5$-geodesic. So, $y_4 - y_5 < \frac{s}{2}.$ If $x_4 > x_6$, then $u_6$ is on a $u_4,u_5$-geodesic, so $x_6>x_4$. If $y_1 \leq y_6$, then $u_6$ is on a $u_1,u_4$-geodesic, so $y_1>y_6$. But if $y_1> y_6$, then $u_6$ is on a $u_1,u_3$-geodesic.

\medskip\noindent
Subcase 1b: $x_3 > x_6$.\\
If $y_4 - y_6 \geq \frac{s}{2}$, then $u_3$ is on a $u_4,u_6$-geodesic. So, $y_4-y_6 < \frac{s}{2}$. But now $u_6$ is on a $u_4,u_5$-geodesic.

\medskip\noindent
{\bf Case 2}: ($y_3 - y_6 < \frac{s}{2}$).\\
Subcase 2a: $x_3 < x_6$.\\
If $x_6>x_4$, then $u_4$ is on a $u_3,u_6$-geodesic. So, $x_4>x_6$. Also, $y_4-y_5 > \frac{s}{2}$, since otherwise $u_6$ is on a $u_4,u_5$-geodesic. If $x_3>x_5$, then $u_3$ is on a $u_4,u_5$-geodesic, so $x_3<x_5$.

We observe that $y_1 < y_6$, for otherwise $u_6$ is on a $u_1,u_5$-geodesic. If $x_6-x_2 > \frac{r}{2}$, then $u_1$ is on a $u_2,u_6$-geodesic. Therefore, $x_6-x_2 < \frac{r}{2}.$ This implies that $y_2>y_5$, for otherwise $u_5$ lies on a $u_2,u_6$-geodesic. But now, $u_5$ lies on a $u_1,u_3$-geodesic.

\medskip\noindent
Subcase 2b ($x_3 > x_6$).\\
If $y_4-y_5 \geq \frac{s}{2}$, then $u_3$ is on a $u_4,u_5$-geodesic, and if  $y_4-y_5 < \frac{s}{2}$, then $u_6$ is on a $u_4,u_5$-geodesic.
\qed

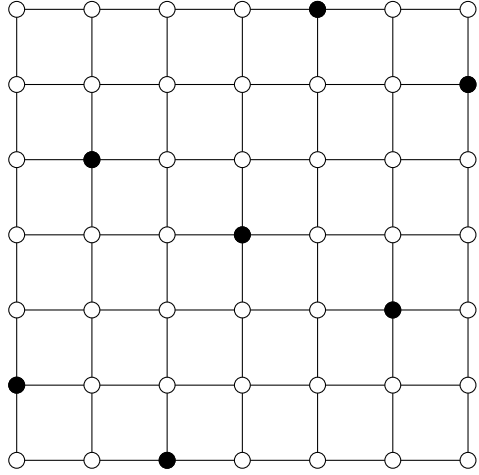
\begin{figure}[ht!]
\begin{center}
\begin{tikzpicture}
\foreach \x in {0,1,...,6}
    {\draw (\x,0) -- (\x,6);}
\foreach \x in {0,1,...,6}
    {\draw (0,\x) -- (6,\x);}
\foreach \x in {0,1,...,6}
\foreach \y in {0,1,...,6}
\filldraw[fill=white] (\x,\y) circle (3pt);
\filldraw[fill=black] (0,1) circle (3pt);
\filldraw[fill=black] (1,4) circle (3pt);
\filldraw[fill=black] (2,0) circle (3pt);
\filldraw[fill=black] (3,3) circle (3pt);
\filldraw[fill=black] (4,6) circle (3pt);
\filldraw[fill=black] (5,2) circle (3pt);
\filldraw[fill=black] (6,5) circle (3pt);
\end{tikzpicture}
\end{center}
    \caption{General position set with seven vertices.}
    \label{fig:seven_vertex}
\end{figure}

A general position set $S$ of $C_7 \cp C_7$ with $|S| = 7$ is shown on Fig.~\ref{fig:seven_vertex}. It can be easily checked that the minimum distance between the pairs of vertices from $S$ is $3$, and that the maximal distance is $5$. Therefore, $S$ is indeed a general position set.

The construction from Fig.~\ref{fig:seven_vertex} cannot be extended to arbitrary $r,s \ge 7$. We have checked by computer that $\gp(C_8\cp C_7)\leq 6$. In fact, $6$ is also the exact lower bound in this case as follows from the following, main result of this section. 

\begin{theorem}\label{thm:torus}
If $r\ge s \ge 3$, $s\ne 4$, and $r\ge 6$, then $\gp(C_r\cp C_s) \in \{6, 7\}$.
\end{theorem}

\proof
From Theorem~\ref{thm:torus_up} we know that $\gp(C_r\cp C_s)\leq 7$. It thus remains to prove that $\gp(C_r\cp C_s)\ge 6$.  

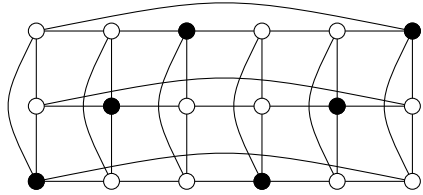
\begin{figure}[ht!]
\begin{center}
\begin{tikzpicture}
\foreach \y in {1,2,3}
{\draw (1,\y)..controls+(2.5,0.5) ..(6,\y);}
\foreach \x in {1,...,6}
{\draw (\x,1)..controls+(-0.5,1) ..(\x,3);}
\foreach \x in {1,...,6}
\foreach \y in {1,...,2}
{\draw (\x,\y) -- (\x,\y+1); }
\foreach \x in {1,...,5}
\foreach \y in {1,...,3}
{\draw (\x,\y) -- (\x+1,\y); }
\foreach \x in {1,...,6}
\foreach \y in {1,...,3}
{
\filldraw[fill=white] (\x,\y) circle (3pt);

}
\foreach \x in {1,2,3}
{\filldraw[black] (\x,\x) circle (3pt);
\filldraw[black] (\x+3,\x) circle (3pt);}
\end{tikzpicture}\caption{The set $S$ in $C_6\cp C_3$ appears in bold.}
\label{fig:C6-times-C3}
\end{center}
\end{figure}

The condition $s\ne 4$ assures that $S_s = \{0,\lfloor\frac{s}{3}\rfloor,\lfloor\frac{2s}{3}\rfloor\}$ is a gp-set of $C_s$. The condition that $r\ge 6$ assures that $\lfloor r/6\rfloor \ge 1$. Consider now the set
\begin{align*}
S = \{ & (0,0),
(\left\lfloor r/2\right\rfloor,0),
(\left\lfloor r/6 \right\rfloor, \left\lfloor s/3 \right\rfloor),
(\left\lfloor r/6 \right\rfloor + \left\lfloor r/2\right\rfloor,\left\lfloor s/3 \right\rfloor), \\
& (\left\lfloor (2r)/6 \right\rfloor, \left\lfloor (2s)/3 \right\rfloor),
(\left\lfloor (2r)/6 \right\rfloor + \left\lfloor r/2\right\rfloor, \left\lfloor (2s)/3 \right\rfloor)\}\,.
\end{align*}
In Fig.~\ref{fig:C6-times-C3} the set $S$ is shown for the case $C_6\cp C_3$.

We claim that $S$ is a general position set. Since $C_p$-layers are isometric subgraphs, no other vertex is on a geodesic between the pair of vertices with the same second coordinate. Hence we only need to consider the triples of vertices from $S$ with pairwise different second coordinates. We do this for the vertices $x_1 = (0,0)$, $x_2 = (\left\lfloor r/6 \right\rfloor, \left\lfloor s/3 \right\rfloor)$, and $x_3 = (\left\lfloor (2r)/6 \right\rfloor, \left\lfloor (2s)/3 \right\rfloor)$, the other cases are treated similarly. Since $x_1, x_2, x_3$ lie in a subgraph of $C_r\cp  C_s$ isomorphic to $P_{\lfloor r/3 \rfloor + 1} \cp C_s$ which is an isometric subgraph, it suffices to show that $d(x_1,x_3) < d(x_1,x_2) + d(x_2,x_3)$. This can be verified using the facts
$d(x_1,x_3) = \lfloor r/3 \rfloor + (s - \lfloor (2s)/3 \rfloor)$,
$d(x_1,x_2) = \lfloor r/6 \rfloor + \lfloor s/3 \rfloor$, and
$d(x_2,x_3) \ge \lfloor r/6 \rfloor + \lfloor s/3 \rfloor$.
\qed

Theorems~\ref{thm:torus_up} and~\ref{thm:torus} yield the following problem, for which we believe that if ${\rm min}\{r,s\} \leq 6$, then $\gp(C_r \cp C_s) = 6$. 

\begin{problem}
Determine for $r\ge 3, r\neq 4, s\ge 6$, whether $\gp(C_r\cp C_s) = 6$ or whether $\gp(C_r\cp C_s) = 6$.
\end{problem}

\section{Cartesian powers}
\label{sec:3D-and-hypercubes}

In this section, we consider the general position number of Cartesian powers and obtain asymptotically exponential lower bounds using a probabilistic approach.

The $n$-dimensional hypercube $Q_n$ is defined as $K_2^{\cp, n}$. In particular, $Q_1 = K_2$, $Q_2 = C_4$, and $Q_3$ is the  graph of the 3-D cube. Cartesian products of complete graphs, known as Hamming graphs, form a natural generalization of hypercubes. In~\cite{ghorbani-2019}, it was proved that if $k\ge 2$ and $n_1, \ldots, n_k\ge 2$, then
\begin{equation}
\label{eq:Hamming}
\gp(K_{n_1}\cp \cdots \cp K_{n_k}) \ge n_1 + \cdots + n_k - k\,.
\end{equation}
Moreover, this lower bound is sharp on products of two complete graphs, that is, $\gp(K_{n_1}\cp K_{n_2}) = n_1 + n_2 - 2$.

The situation above changes dramatically as $k$ grows. K\"orner \cite{korner-1995} obtained a probabilistic construction of general position sets in $Q_n$ of size $\frac{1}{2}\frac{2^n}{\sqrt{3^n}}$. He also pointed out that the problem of finding the size of the largest point set in general position in $Q_n$ is equivalent to finding the largest size of what is called a $(2,1)$-separating system in coding theory. (For more on separating systems, see~\cite{cohen-2003}.) K\"orner was interested in
$$\alpha=\limsup_{n\rightarrow \infty}\frac{\log_2\gp(Q_n)}{n}\,.$$
His probabilistic lower bound gives $\alpha\ge 1-\frac{1}{2}\log_23$ and he also proved $\alpha\le 1/2$. Later, Randriambololona \cite{randriambololona-2013} improved the lower bound to $\alpha\ge \frac{3}{50}\log_211$ with an explicit construction.

The first moment method can be applied in a general setting to obtain large general position sets. For any graph $G$, let $p(G)$ denote the probability that if one picks a triple $(x,y,z)\in V(G)^3$ uniformly at random, then $d_G(y,z)=d_G(y,x)+d_G(x,z)$ holds. Let us call such triples \textit{bad}. Note that this is never the case if $x\neq y$ and $y=z$, so $p(G)\le 1-\frac{|V(G)|-1}{|V(G)|^2}<1$. Let $H=H_1\cp \cdots \cp H_k$. Observe that  the triple $\mathbf{x}=(x_1,\dots,x_k),\mathbf{y}=(y_1,\dots,y_k),\mathbf{z}=(z_1,\dots,z_k)\in V(H)$ is bad in $H$ if and only if the triples $x_i,y_i,z_i$ are bad in $H_i$ for all $i \in [k]$. So if we pick $M$ vertices uniformly at random with repetition from $V(H)$, then the expected value $\mathbb{E}(X)$ of the number $X=X(M)$ of unordered triples on a geodesic will be $3\binom{M}{3}\prod_{i=1}^kp(H_i)$. If $X\le M/2$, then removing one vertex from every bad triple will leave us a general position set of size at least $M/2$. As there is always an instance for which $X\le \mathbb{E}(X)$ holds, we obtain a general position set of size $M/2$ provided $3\binom{M}{3}\prod_{i=1}^kp(H_i)\le M/2$ holds. Therefore, it seems to be interesting to examine
$$\gp_{\cp}(G):=\limsup_{n\rightarrow \infty}\frac{\log_{|V(G)|} \gp(G^{\cp, n})}{n}\,.$$ Clearly, we have $\gp_{\cp}(G)\le 1$ and the above reasoning yields the following theorem.

For a graph $G$ one can consider its Cartesian power $G^{\cp, n}$. Then the required inequality is $3\binom{M}{3}p(G)^n\le M/2$ which is equivalent to $(M-1)(M-2)\le p(G)^{-n}$. Thus there exists a general position set in $G^{\cp, n}$ of size $\frac{1}{2}p(G)^{-n/2}$.
This and the inequality $p(G)\le 1-\frac{|V(G)|-1}{|V(G)|^2}$ yields the following statement.

\begin{theorem}\label{random}
If $G$ is a graph, then
$$\gp_{\cp} (G)\ge \log_{|V(G)|}p(G)^{-1/2}\ge 1-\log_{|V(G)|}(|V(G)|^2-|V(G)|+1)\,.$$
\end{theorem}

Let us calculate $p(G)$ for some graphs. First of all, $p(K_n)=\frac{2n-1}{n^2}$ as in $K_n$ the equality $d(y,z)=d(y,x)+d(x,z)$ holds if and only if $x=y$ or $x=z$. (The case $p(K_2)=\frac{3}{4}$ in Theorem ~\ref{random} is just K\"orner's result.) For even cycles we have $p(C_{2k})=\frac{k(k+3)-1}{4k^2}$. If the vertices are $\{-(k-1),-(k-2),\dots,0,,\dots,k-1,k\}$ in this cyclic order, then by symmetry we can assume $x=0$. There are $4k-1$ triples with $x=y$ or $x=z$ that form bad triples. If $y=k$ or $z=k$, then there are no other bad triples, otherwise for any $y$, there are $k-|y|$ ways to choose $z$ to obtain a bad triple. Similarly, one can verify $p(C_{2k+1})=\frac{k(k+3)+1}{(2k+1)^2}$. Finally, consider the star $S_k$ with $k$ leaves. Then conditioning on whether $x$ is the center or not one obtains $p(S_k)=\frac{1}{k+1}+\frac{k}{k+1}\frac{2k+1}{(k+1)^2}$. Observe that if one picks uniformly at random only among the leaves of $S_k$, then the probability of picking a bad triple is $p'(S_k)=\frac{2k-1}{k^2}$ which for large enough $k$s is roughly $2/3$ of $p(S_k)$, so in this way one obtains the better bound $\gp_{\cp}(S_k)\ge \log_2p'(S_k)^{-1/2}$.

Concerning $\gp_{\cp}(G)$ we wonder whether one can write limit instead of limit superior in the definition of $\gp_{\cp}(G)$. Moreover, by the above we have $\lim_{k\rightarrow\infty}p(C_k)=\frac{1}{4}$. We also pose:

\begin{problem}
Decide whether $\liminf_{k\rightarrow \infty}\frac{\gp_{\cp}(C_k)}{\log_k2}>1$ holds.
\end{problem}

\section*{Acknowledgements}

We acknowledge the financial support from the Slovenian Research Agency (research core funding No.\ P1-0297 and projects J1-9109, J1-1693, N1-0095, N1-0108). This research was done while the last author was visiting the University of Ljubljana, Slovenia, supported by ``Ministerio de Educaci\'on, Cultura y Deporte'', Spain, under the ``Jos\'e Castillejo'' program for young researchers (reference number: CAS18/00030). The research of the second author was supported by the National Research, Development and Innovation Office - NKFIH under the grants SNN 129364 and FK 132060.


\end{document}